\documentclass[12pt]{article}
\usepackage{amsmath}
\usepackage{amsfonts}
\newtheorem{theorem}{Theorem}
\newtheorem{lemma}{Lemma}
\def\qed{\quad \rule[-.5ex]{1.5ex}{1.5ex} \bigskip}

\begin{document}
\title{Some classical multiple orthogonal polynomials \thanks{The first author is a Research Director of the Belgian Fund for Scientific Research (FWO). Research supported by INTAS and FWO research project 
G.0278.97.}}
\author{Walter Van Assche and Els Coussement
\\
\normalsize{Department of Mathematics, Katholieke Universiteit Leuven}}
\date{}
\maketitle

\section{Classical orthogonal polynomials}
One aspect in the theory of orthogonal polynomials is their study as
special functions. Most important orthogonal polynomials can be written
as terminating hypergeometric series and during the twentieth century
people have been working on a classification of all such hypergeometric
orthogonal polynomial and their characterizations.

The \textit{very classical orthogonal polynomials} are those named after
Jacobi, Laguerre, and Hermite. In this paper we will always be considering
monic polynomials, but in the literature one often uses a different
normalization. Jacobi polynomials are (monic) polynomials of degree $n$ which are orthogonal to all lower degree polynomials with respect to the weight function $(1-x)^\alpha (1+x)^\beta$ on $[-1,1]$, where
$\alpha, \beta > -1$. The change of variables $x \mapsto 2x-1$ gives
Jacobi polynomials on $[0,1]$ for the weight function $w(x) = x^\beta(1-x)^\alpha$, and we will denote these (monic) polynomials by 
$P_n^{(\alpha,\beta)}(x)$. They are defined by the orthogonality conditions
\begin{equation}  \label{eq:jacobi}
  \int_0^1 P_n^{(\alpha,\beta)}(x) x^\beta(1-x)^\alpha x^k\, dx = 0,
   \qquad k=0,1,\ldots, n-1. 
\end{equation}
The monic Laguerre polynomials $L_n^{(\alpha)}(x)$ (with $\alpha > -1$)
are orthogonal on $[0,\infty)$ to all polynomials of degree less than $n$
with respect to the weight $w(x) =x^\alpha e^{-x}$ and hence satisfy
the orthogonality conditions
\begin{equation}  \label{eq:laguerre}
  \int_0^\infty L_n^{(\alpha)}(x) x^\alpha e^{-x} x^k\, dx = 0, 
 \qquad k=0,1,\ldots,n-1.
\end{equation}
Finally, the (monic) Hermite polynomials $H_n(x)$ are orthogonal to all lower degree polynomials with respect to the weight function
$w(x) = e^{-x^2}$ on $(-\infty,\infty)$, so that
\begin{equation}
  \int_{-\infty}^\infty H_n(x) e^{-x^2} x^k\, dx = 0, \qquad
  k=0,1,\ldots,n-1.
\end{equation}
These three families of orthogonal polynomials can be characterized
in a number of ways:
\begin{itemize}
\item Their weight functions $w$ satisfy a first order differential
equation with polynomial coefficients
\begin{equation} \label{eq:pearson}
   \sigma(x) w'(x) = \rho(x) w(x),
\end{equation}
with $\sigma$ of degree at most two and $\rho$ of degree one. This
equation is known as \textit{Pearson's equation} and also appears
in probability theory, where the corresponding weights (densities) are
known as the beta density (Jacobi), the gamma density (Laguerre), and
the normal density (Hermite). Note however that for probability
density functions one needs to normalize these weights appropriately.
For the Jacobi weight we have $\sigma(x) = x(1-x)$, for the Laguerre
weight we have $\sigma(x) = x$, and for the Hermite weight we see that
$\sigma(x) = 1$, so that each family corresponds to a different degree
of the polynomial $\sigma$.
\item The derivatives of the very classical polynomials are again
orthogonal polynomials of the same family but with different parameters
(Sonin 1887, W. Hahn 1949).
Indeed, integration by parts of the orthogonality relations and the use
of Pearson's equation show that
\begin{eqnarray*}
   \frac{d}{dx} P_n^{(\alpha,\beta)}(x) & = & 
  n P_{n-1}^{(\alpha+1,\beta+1)}(x), \\
   \frac{d}{dx} L_n^{(\alpha)}(x) & = & n L_{n-1}^{(\alpha+1)}(x), \\
  \frac{d}{dx} H_n(x) & = & n H_{n-1}(x).
\end{eqnarray*}
The differential operator $D = d/dx$ therefore acts as a \textit{lowering
operator} that lowers the degree of the polynomial.
\item Pearson's equation also gives rise to a raising operator that
raises the degree of the polynomials. Indeed, integration by parts
shows that
\begin{eqnarray}
  \frac{d}{dx} \left[x^\beta(1-x)^\alpha P_n^{(\alpha,\beta)}(x)\right]
 & = & -(\alpha+\beta+n) x^{\beta-1}(1-x)^{\alpha-1}
  P_{n+1}^{(\alpha-1,\beta-1)}(x), \\
 \frac{d}{dx} \left[ x^\alpha e^{-x} L_n^{(\alpha)}(x) \right] & = & 
  - x^{\alpha-1} e^{-x} L_{n+1}^{(\alpha-1)}(x), \\
  \frac{d}{dx} \left[ e^{-x^2} H_n(x) \right] & = & -2 e^{-x^2} H_{n+1}(x).
\end{eqnarray}
The raising operator is therefore of the form $\sigma(x)/w(x) D w(x)$.
Using this raising operation repeatedly gives the \textit{Rodrigues
formula} for these orthogonal polynomials: 
\begin{eqnarray}
 \frac{d^n}{dx^n} \left[x^{\beta+n}(1-x)^{\alpha+n} \right] & = & (-1)^n
 (\alpha+\beta+n+1)_n x^\beta(1-x)^\alpha
  P_n^{(\alpha,\beta)}(x), \label{eq:jacobiRod} \\ 
  \frac{d^n}{dx^n} \left[ x^{\alpha+n} e^{-x} \right]
  & = & (-1)^n x^\alpha e^{-x} L_n^{(\alpha)}(x), \\
  \frac{d^n}{dx^n} e^{-x^2} & = & (-1)^n 2^n e^{-x^2} H_n(x) .
\end{eqnarray}
The Rodrigues formula is therefore of the form
\[  \frac{d^n}{dx^n}\left[ \sigma^n(x)w(x) \right] = C_n w(x) P_n(x), \]
where $C_n$ is a normalization constant (Hildebrandt 1931).
\item Combining the lowering and the raising operator gives a linear second
order differential equation for these orthogonal polynomials, of the form
\begin{equation}  \label{eq:diffeq}
   \sigma(x) y''(x) + \tau(x) y'(x) = \lambda_n y(x), 
\end{equation}
where $\sigma$ is a polynomial of degree at most 2 and $\tau$
a polynomial of degree at most 1, both independent of the degree $n$, and
$\lambda_n$ is a constant depending on $n$ (Bochner 1929).
\end{itemize}

The Laguerre polynomials and the Hermite polynomials are limiting
cases of the Jacobi polynomials. Indeed, one has
\begin{equation}  \label{eq:jactolag}
 \lim_{\alpha \to \infty} \alpha^n P_n^{(\alpha,\beta)}(x/\alpha)
   = L_n^{(\beta)} (x), 
\end{equation}
and
\begin{equation}  \label{eq:jactoher}
  \lim_{\alpha \to \infty}  2^n \alpha^{n/2} P_n^{(\alpha,\alpha)}(\frac{x+\sqrt{\alpha}}{2\sqrt{\alpha}})= H_n(x). 
\end{equation}
The Hermite polynomials are also a limit case of the Laguerre
polynomials:
\begin{equation}  \label{eq:lagtoher}
   \lim_{\alpha \to \infty} (2\alpha)^{-n/2} L_n^{(\alpha)}(\sqrt{2\alpha}x+\alpha) = H_n(x).
\end{equation}
In this respect the Jacobi, Laguerre and Hermite polynomials are in a hierarchy, with Jacobi leading to Laguerre and Laguerre leading to Hermite,
and with a shortcut for Jacobi leading to Hermite. This is just a
very small piece in a large table known as Askey's table which also
contains classical orthogonal polynomials of a discrete variable
(Hahn, Meixner, Kravchuk, and Charlier) for which 
the differential operator $D$ needs to be replaced by difference
operators $\Delta$ and $\nabla$ on a linear lattice (a lattice
with constant mesh, see \cite{nikiuvasus}). Finally,
allowing a quadratic lattice also gives Meixner-Pollaczek, dual Hahn,
continuous Hahn, continuous dual Hahn, Racah, and Wilson polynomials,
which are all in the Askey table. These polynomials have a number of 
$q$-extensions involving the $q$-difference operator and leading to the
$q$-extension of the Askey table. 
In \cite{andask} Andrews and Askey suggest to define the classical
orthogonal polynomials as those polynomials that are a
limiting case of the ${}_4\varphi_3$-polynomials
\[  R_n(\lambda(x);a,b,c,d;q) =
 {}_4\varphi_3 \left(   \begin{array}{c}
    q^{-n}, q^{n+1} ab, q^{-x}, q^{x+1} cd \\
    aq, bdq, cq \end{array} ; q,q  \right), \]
with $\lambda(x) = q^{-x}+q^{x+1}cd$ and $bdq = q^{-N}$ (these are
the $q$-Racah polynomials) or the ${}_4\varphi_3$-polynomials
\[ \frac{a^n W_n(x;a,b,c,d|q)}{(ab;q)_n (ac;q)_n(ad;q)_n} =
  {}_4\varphi_3 \left(   \begin{array}{c}
    q^{-n}, q^{n-1} abcd, ae^{i\theta}, ae^{-i\theta} \\
    ab, ac, ad \end{array} ; q,q  \right), \]
with $x=\cos \theta$ (these are the Askey-Wilson polynomials).
All these classical orthogonal polynomials then have the following properties:
\begin{itemize}
\item they have a Rodrigues formula,
\item an appropriate divided difference operator acting on them gives a set
of orthogonal polynomials,
\item they satisfy a second order difference equation in $x$ which is
of Sturm-Liouville type.
\end{itemize}
The classical orthogonal polynomials in this wide sense have been the subject of intensive research during 
the twentieth century.
We recommend the report by Koekoek and Swarttouw \cite{koekswart},
the book by Andrews, Askey and Roy \cite{andaskroy}, and the books
by Nikiforov,  Uvarov \cite{nikiuva}, and 
Nikiforov, Suslov, and Uvarov \cite{nikiuvasus}
for more material. Szeg\H{o}'s book \cite{szego} is still a very good source for the very classical orthogonal polynomials of Jacobi, Laguerre, and Hermite. For characterization results one should consult a survey
by Al-Salam \cite{alsalam}.

\section{Multiple orthogonal polynomials}
Recently, there has been a renewed interest in an extension of the notion 
of orthogonal polynomials known as multiple orthogonal polynomials. This
notion comes from simultaneous rational approximation, in particular from
Hermite-Pad\'e approximation of a system of $r$ functions, and hence
has its roots in the nineteenth century. 
However, only recently examples of multiple orthogonal polynomials
appeared in the (mostly Eastern European) literature.
In this paper
we will introduce multiple orthogonal polynomials using the orthogonality
relations and we will only use weight functions. The extension to
measures is straightforward.

Suppose we are given $r$ weight functions $w_1,w_2,\ldots,w_r$ on the real
line and that the support of each $w_i$ is a subset of an interval
$\Delta_i$. We will often be using a multi-index $\vec{n} = (n_1,n_2,\ldots,n_r) \in \mathbb{N}^r$ and its length
$|\vec{n}| = n_1+n_2+\cdots+n_r$. 
\begin{itemize}
\item
The $r$-vector of \textit{type I multiple orthogonal polynomials} $(A_{\vec{n},1},\ldots,A_{\vec{n},r})$ is such that each $A_{\vec{n},i}$
is a polynomial of degree $n_i-1$ and the following orthogonality
conditions hold:
\begin{equation} \label{eq:Iorth}
  \int x^k \sum_{j=1}^r A_{\vec{n},j}(x) w_j(x) \, dx = 0, \qquad
  k=0,1,2,\ldots,|\vec{n}|-2. 
\end{equation}
Each $A_{\vec{n},i}$ has $n_i$ coefficients so that the type I vector
is completely determined if we can find all the $|\vec{n}|$ unknown
coefficients. The orthogonality relations (\ref{eq:Iorth}) give
$|\vec{n}|-1$ linear and homogeneous relations for these $|\vec{n}|$
coefficients. If the matrix of coefficients has full rank, then
we can determine the type I vector uniquely up to a multiplicative factor.
\item
The type II multiple orthogonal polynomial $P_{\vec{n}}$ is the polynomial
of degree $|\vec{n}|$ that satisfies the following orthogonality
conditions
\begin{eqnarray}
 \int_{\Delta_1} P_{\vec{n}}(x) w_1(x) x^k\, dx & = & 0, 
   \qquad k=0,1,\ldots,n_1-1,  \label{eq:IIorth1} \\
 \int_{\Delta_2} P_{\vec{n}}(x) w_2(x) x^k\, dx & = & 0, 
   \qquad k=0,1,\ldots,n_2-1,  \label{eq:IIorth2} \\
   & \vdots &  \nonumber \\
 \int_{\Delta_r} P_{\vec{n}}(x) w_r(x) x^k\, dx & = & 0, 
   \qquad k=0,1,\ldots,n_r-1. \label{eq:IIorthr}
\end{eqnarray}
This gives $|\vec{n}|$ linear and homogeneous equations for the
$|\vec{n}|+1$ unknown coefficients of $P_{\vec{n}}(x)$. We will choose
the type II multiple orthogonal polynomials to be monic so that the
remaining $|\vec{n}|$ coefficients can be determined uniquely by the
orthogonality relations, provided the matrix of coefficients has full rank.
\end{itemize}

In this paper the emphasis will be on type II multiple orthogonal
polynomials. The unicity of multiple orthogonal polynomials can
only be guaranteed under additional assumptions on the $r$ weights.
Two distinct cases for which the type II multiple
orthogonal polynomials are given as follows.
\begin{enumerate}
\item In an \textit{Angelesco system} (Angelesco, 1918) the intervals $\Delta_i$, on which
the weights are supported, are disjoint, i.e., $\Delta_i \cap \Delta_j
= \emptyset$ whenever $i \neq j$. Actually, it is sufficient that the
open intervals $\stackrel{\circ}{\Delta}_i$ are disjoint, so that
the closed intervals $\Delta_i$ are allowed to touch. 

\begin{theorem}
In an Angelesco
system the type II multiple orthogonal polynomial $P_{\vec{n}}(x)$
factors into $r$ polynomials $\prod_{j=1}^r q_{n_j}(x)$, where
each $q_{n_j}$ has exactly $n_j$ zeros on $\Delta_j$. 
\end{theorem}
\textbf{Proof:}
Suppose $P_{\vec{n}}(x)$ has $m_j < n_j$ sign changes on $\Delta_j$
at the points $x_1,\ldots,x_{m_j}$. Let $Q_{m_j}(x) = (x-x_1)\cdots
(x-x_{m_j})$, then $P_{\vec{n}}(x)Q_{m_j}(x)$ does not change
sign on $\Delta_j$, and hence
\[   \int_{\Delta_j} P_{\vec{n}}(x)Q_{m_j}(x) w_j(x)\, dx \neq 0. \]
But this is in contradiction with the orthogonality relation on
$\Delta_j$. Hence $P_{\vec{n}}(x)$ has at least $n_j$ zeros on
$\Delta_j$. Now all the intervals $\Delta_j$ $(j=1,2,\ldots,r)$
are disjoint, hence this gives at least $|\vec{n}|$ zeros of
$P_{\vec{n}}(x)$ on the real line. The degree of this polynomials is
precisely $|\vec{n}|$, so there are exactly $n_j$ zeros on each interval
$\Delta_j$. \qed

\item For an \textit{AT system} all the weights are supported on the
same interval $\Delta$, but we require that the $|\vec{n}|$ functions 
\begin{multline*}
 w_1(x), x w_1(x), \ldots, x^{n_1-1}w_1(x), w_2(x), xw_2(x), \ldots,
x^{n_2-1}w_2(x),  \\
  \ldots, w_r(x), xw_r(x),\ldots,x^{n_r-1}w_r(x)
\end{multline*}
form
a Chebyshev system on $\Delta$ for each multi-index $\vec{n}$. This means
that every linear combination 
\[   \sum_{j=1}^r Q_{n_j-1}(x) w_j(x), \]
with $Q_{n_j-1}$ a polynomial of degree at most $n_j-1$, has at most
$|\vec{n}|-1$ zeros on $\Delta$.

\begin{theorem}
In an AT system the type II multiple orthogonal polynomial $P_{\vec{n}}(x)$
has exactly $|\vec{n}|$ zeros on $\Delta$. For the type I vector of multiple orthogonal polynomials, the linear combination
$\sum_{j=1}^r A_{\vec{n},j}(x)w_j(x)$ has exactly $|\vec{n}|-1$ zeros on $\Delta$.
\end{theorem}
\textbf{Proof:}
Suppose $P_{\vec{n}}(x)$ has $m < |\vec{n}|$ sign changes on
$\Delta$ at the points $x_1,\ldots,x_m$. Take a multi-index
$\vec{m}=(m_1,m_2,\ldots,m_r)$ such that $m_i\leq n_i$ for every
$i$ and $m_j < n_j$ for some $j$ and construct the function
\[  Q(x) = \sum_{i=0}^r Q_i(x)w_i(x) , \]
where each $Q_i$ is a polynomial of degree $m_i-1$ whenever $i \neq j$, and
$Q_j$ is a polynomial of degree $m_j$, satisfying the interpolation conditions
\[    Q(x_k) = 0, \qquad k=1,2,\ldots,m, \]
and $Q(x_0)=1$ for an additional point $x_0 \in \Delta$. This interpolation
problem has a unique solution since we are dealing with a Chebyshev
system. The function $Q$ has already $m$ zeros, and since we are in
a Chebyshev system, it can have no additional sign changes. Furthermore,
the function does not vanish identically since $Q(x_0)=1$. Obviously
$P_{\vec{n}}(x)Q(x)$ does not change sign on $\Delta$, so that
\[  \int_\Delta P_{\vec{n}}(x)Q(x)\, dx \neq 0, \]
but this is in contrast with the orthogonality relations for the
type II multiple orthogonal polynomial. Hence $P_{\vec{n}}(x)$ has
exactly $|\vec{n}|$ zeros on $\Delta$.

The proof for the type I multiple orthogonal polynomials is similar.
First of all, since we are dealing with an AT system, the function
\[  A(x) = \sum_{j=1}^r A_{\vec{n},j}(x)w_j(x) \]
has at most $|\vec{n}|-1$ zeros on $\Delta$. Suppose it has
$m < |\vec{n}|-1$ sign changes at the points $x_1,x_2,\ldots,x_m$, then
we use the polynomial $Q_m(x)=(x-x_1)\cdots(x-x_m)$ so that
$A(x)Q(x)$ does not change sign on $\Delta$, and
\[  \int_\Delta  A(x) Q(x) \, dx \neq 0, \]
which is in contradiction with the orthogonality of the type I multiple
orthogonal polynomial. Hence $A(x)$ has exactly $|\vec{n}|-1$ zeros
on $\Delta$. \qed

\end{enumerate}

Orthogonal polynomials on the real line always satisfy a three-term
recurrence relation. There are also finite order recurrences for
multiple orthogonal polynomials, and there are quite a few of recurrence relations possible since we are dealing with multi-indices. There is
an interesting recurrence relation of order $r+1$ for the type II multiple orthogonal polynomials with nearly diagonal 
multi-indices. Let $n \in \mathbb{N}$ and write it as
$n = kr+j$, with $0 \leq j < r$. The nearly diagonal multi-index $\vec{s}(n)$
corresponding to $n$ is then given by 
\[  \vec{s}(n) = (\underbrace{k+1,k+1,\ldots,k+1}_{j\ \mathrm{times}},
\underbrace{k,k,\ldots,k}_{r-j\ \mathrm{times}}). \]
If we denote the corresponding multiple orthogonal polynomials by
\[  P_n(x) = P_{\vec{s}(n)}(x), \]
then the following recurrence relation holds:
\begin{equation}
  xP_n(x) = P_{n+1}(x) + \sum_{j=0}^r a_{n,j} P_{n-j}(x), 
\end{equation}
with initial conditions $P_0(x) = 1$, $P_j(x)=0$ for $j=-1,-2,\ldots,-r$.
The matrix
\[  \begin{pmatrix} 
  a_{0,0} & 1 &    \\
  a_{1,1} & a_{1,0} & 1 \\
  a_{2,2} & a_{2,1} & a_{2,0} & 1 \\
   \vdots & & & \ddots &\ddots \\
  a_{r,r} & a_{r,r-1} & \cdots & & a_{r,0} & 1 \\
          & a_{r+1,r} & \ddots &  & & a_{r+1,0} & 1 \\
          &   & \ddots & \ddots & &  & \ddots & \ddots\\
           &  &   &  \ddots & \ddots & & & \ddots & 1 \\
          &   &   &     & a_{n,r} & a_{n,r-1} & \cdots & a_{n,1} &  a_{n,0}
 \end{pmatrix} \]
has eigenvalues at the zeros of $P_{n+1}(x)$, so that in the case of
Angelesco systems or AT systems we are dealing with non-symmetric matrices
with real eigenvalues. The infinite matrix will act as an operator
on $\ell^2$, but this operator is never self-adjoint and furthermore
has not a simple spectrum, as is the case for ordinary orthogonal
polynomials. Now there will be a set of $r$ cyclic vectors and the
spectral theory of this operator becomes more complicated (and more
interesting). There are many open problems concerning this non-symmetric
operator.

\section{Some very classical multiple orthogonal polynomials}
We will now describe seven families of multiple orthogonal polynomials
which have the same flavor as the very classical orthogonal
polynomials of Jacobi, Laguerre, and Hermite.
They certainly deserve to be called classical since they have
a Rodrigues formula and there is a first order differential operator which, when applied to these classical multiple orthogonal polynomials,
gives another set of multiple orthogonal polynomials. However, these
are certainly not the only families of multiple orthogonal polynomials
(see Section 4.1).  The first four families
are AT systems which are connected by limit passages, the last three
families are Angelesco systems which are also connected by limit
passages. All these families have been introduced in the literature before.
We will list some of their properties and give explicit formulas, most
of which have not appeared earlier.

\setlength{\unitlength}{4pt}
\begin{picture}(120,90)
\put(21,75){AT systems}
\put(75,75){Angelesco systems}
\put(17,60){\framebox(24,10){\shortstack{Jacobi-Pi\~neiro \\ $P_{n,m}^{(\alpha_0,\alpha_1,\alpha_2)}(x)$}}}
\put(75,60){\framebox(24,10){\shortstack{Jacobi-Angelesco \\ $P_{n,m}^{(\alpha,\beta,\gamma)}(x;a)$}}}
\put(0,40){\framebox(28,10){\shortstack{multiple Laguerre I \\
$L_{n,m}^{(\alpha_1,\alpha_2)}(x)$}}}
\put(30,40){\framebox(27,10){\shortstack{multiple Laguerre II \\
$L_{n,m}^{(\alpha_0,c_1,c_2)}(x)$}}}
\put(17,20){\framebox(24,10){\shortstack{multiple Hermite \\
$H_{n,m}^{(c_1,c_2)}(x)$}}}
\put(60,40){\framebox(25,10){\shortstack{Jacobi-Laguerre \\ 
$L_{n,m}^{(\alpha,\beta)}(x;a)$}}}
\put(90,40){\framebox(25,10){\shortstack{Laguerre-Hermite \\
$H_{n,m}^{(\beta)}(x)$}}}
\put(29,60){\vector(0,-1){28}}
\put(29,60){\vector(-2,-1){15}}
\put(29,60){\vector(2,-1){15}}
\put(87,60){\vector(-2,-1){15}}
\put(87,60){\vector(2,-1){15}}
\put(14,40){\vector(1,-1){8}}
\end{picture}

\subsection{Jacobi-Pi\~neiro polynomials}
The Jacobi-Pi\~neiro polynomials are multiple orthogonal polynomials
associated with an AT system consisting of Jacobi weights on $[0,1]$
with different singularities at $0$ and the same singularity at $1$.
They were first studied by Pi\~neiro \cite{pineiro} when $\alpha_0=0$.
The general case appears in \cite[p.~162]{nikisor}.
Let $\alpha_0 > -1$ and $\alpha_1,\ldots,\alpha_r$ be such that each
$\alpha_i > -1$ and $\alpha_i - \alpha_j \notin \mathbb{Z}$ whenever
$i \neq j$. The Jacobi-Pi\~neiro polynomial $P_{\vec{n}}^{(\alpha_0,\vec{\alpha})}$ for the multi-index
$\vec{n} = (n_1,n_2,\ldots,n_r) \in \mathbb{N}^r$ and $\vec{\alpha}
= (\alpha_1,\ldots,\alpha_r)$ is the monic polynomial of degree
$|\vec{n}| = n_1+n_2+\cdots+n_r$ that satisfies the orthogonality 
conditions
\begin{eqnarray}
  \int_0^1 P_{\vec{n}}^{(\alpha_0,\vec{\alpha})}(x) x^{\alpha_1}
    (1-x)^{\alpha_0} x^k \, dx & = & 0, \qquad k=0,1,\ldots,n_1-1,
  \label{eq:pineiro1} \\
 \int_0^1 P_{\vec{n}}^{(\alpha_0,\vec{\alpha})}(x) x^{\alpha_2}
    (1-x)^{\alpha_0} x^k \, dx & = & 0, \qquad k=0,1,\ldots,n_2-1, 
  \label{eq:pineiro2} \\
  & \vdots &  \nonumber \\
 \int_0^1 P_{\vec{n}}^{(\alpha_0,\vec{\alpha})}(x) x^{\alpha_r}
    (1-x)^{\alpha_0} x^k \, dx & = & 0, \qquad k=0,1,\ldots,n_r-1.
  \label{eq:pineiror}
\end{eqnarray}
Since each weight $w_i(x) = x^{\alpha_i}(1-x)^{\alpha_0}$ satisfies
a Pearson equation
\[   x(1-x) w_i'(x) = [\alpha_i(1-x)-\alpha_0 x] w_i(x) \]
and the weights are related by
\[   w_i(x) = x^{\alpha_i-\alpha_j} w_j(x), \]
one can use integration by parts on each of the $r$ integrals 
(\ref{eq:pineiro1})--(\ref{eq:pineiror}) to find the following 
raising operators:
\begin{equation}  \label{eq:pineirorais}
   \frac{d}{dx} \left( x^{\alpha_j}(1-x)^{\alpha_0} 
  P_{\vec{n}}^{(\alpha_0,\vec{\alpha})}(x) \right) =
  -( |\vec{n}|+\alpha_0+\alpha_j) x^{\alpha_j-1} (1-x)^{\alpha_0-1}
 P_{\vec{n}+\vec{e}_j}^{(\alpha_0-1,\vec{\alpha}-\vec{e}_j)}(x),
\end{equation}
where $\vec{e}_j$ is the $j$th standard unit vector. Repeatedly using this
raising operator gives the Rodrigues formula
\begin{equation}  \label{eq:pineiroRod}
 (-1)^{|\vec{n}|} \prod_{j=1}^r (|\vec{n}|+\alpha_0+\alpha_j+1)_{n_j}
  P_{\vec{n}}^{(\alpha_0,\vec{\alpha})}(x)
 = (1-x)^{-\alpha_0} \prod_{j=1}^r \left[ x^{-\alpha_j} 
  \frac{d^{n_j}}{dx^{n_j}} x^{n_j+\alpha_j} \right] (1-x)^{\alpha_0+|\vec{n}|}.
\end{equation}
The product of the $r$ differential operators $x^{-\alpha_j} 
  D^{n_j} x^{n_j+\alpha_j}$ on the right hand side
can be taken in any order since these operators are commuting.

The Rodrigues formula allows us to obtain an explicit expression. For the
case $r=2$ we write
\begin{multline}
 (-1)^{n+m} (n+m+\alpha_0+\alpha_1+1)_n  (n+m+\alpha_0+\alpha_2+1)_m P_{n,m}^{(\alpha_0,\alpha_1,\alpha_2)}(x) \\
  = (1-x)^{-\alpha_0} x^{-\alpha_1} \frac{d^n}{dx^n} x^{\alpha_1-\alpha_2+n}
  \frac{d^m}{dx^m} x^{\alpha_2+m} (1-x)^{\alpha_0+n+m}. 
\end{multline}
The $m$th derivative can be worked out using the Rodrigues
formula (\ref{eq:jacobiRod}) for Jacobi polynomials and gives
\begin{multline*}
 (-1)^{n} (n+m+\alpha_0+\alpha_1+1)_n   P_{n,m}^{(\alpha_0,\alpha_1,\alpha_2)}(x) \\
  = (1-x)^{-\alpha_0} x^{-\alpha_1} \frac{d^n}{dx^n} x^{\alpha_1+n}
   (1-x)^{\alpha_0+n}   P_m^{(\alpha_0+n,\alpha_2)}(x). 
\end{multline*}
Now use Leibniz' rule to work out the $n$th derivative:
\begin{multline*}
  (-1)^{n} (n+m+\alpha_0+\alpha_1+1)_n   P_{n,m}^{(\alpha_0,\alpha_1,\alpha_2)}(x) \\
  = (1-x)^{-\alpha_0} x^{-\alpha_1} \sum_{k=0}^n \binom{n}{k}
  \frac{d^k}{dx^k} x^{\alpha_1+n} \frac{d^{n-k}}{dx^{n-k}}
  (1-x)^{\alpha_0+n} P_m^{(\alpha_0+n,\alpha_2)}(x) . 
\end{multline*}
In order to work out the derivative involving the Jacobi polynomial, we will use the following lemma.

\begin{lemma}
Let $P_n^{(\alpha,\beta)}(x)$ be the $n$th degree monic Jacobi
polynomial on $[0,1]$. Then for $\alpha >0$ and $\beta > -1$
\begin{equation}  \label{eq:lem1}
   \frac{d}{dx} \left[ (1-x)^\alpha P_n^{(\alpha,\beta)}(x) \right]
   = -( \alpha+n) (1-x)^{\alpha-1} P_n^{(\alpha-1,\beta+1)}(x),
\end{equation}
and
\begin{equation}  \label{eq:lemm}
   \frac{d^m}{dx^m} \left[ (1-x)^\alpha P_n^{(\alpha,\beta)}(x) \right]
   = (-1)^m ( \alpha+n-m+1)_m (1-x)^{\alpha-m} P_n^{(\alpha-m,\beta+m)}(x).
\end{equation}
\end{lemma}
\textbf{Proof:}
First of all, observe that
\[  \frac{d}{dx} \left[ (1-x)^\alpha P_n^{(\alpha,\beta)}(x) \right]
     = (1-x)^{\alpha-1} \left( -\alpha P_n^{(\alpha,\beta)}(x)
    + (1-x) [P_n^{(\alpha,\beta)}(x)]' \right), \]
so that the right hand side is $-(\alpha+n)(1-x)^{\alpha-1} Q_n(x)$, with
$Q_n$ a monic polynomial of degree $n$. Integrating by parts gives
\begin{multline*}
 -(\alpha+n) \int_0^1 (1-x)^{\alpha-1} x^{\beta+k+1} Q_n(x)\, dx \\
   = \left. x^{\beta+k+1}(1-x)^{\alpha} P_n^{(\alpha,\beta)}(x) \right|_0^1
   - (\beta+k+1) \int_0^1 x^{\beta+k}(1-x)^{\alpha} P_n^{(\alpha,\beta)}(x)\, dx . 
\end{multline*}
Obviously, when $\alpha > 0$ and $\beta > -1$, then the integrated terms
on the right hand side vanish. The integral on the right hand side
vanishes for $k=0,1,\ldots,n-1$ because of orthogonality. Hence
$Q_n$ is a monic polynomial which is orthogonal to all polynomials
of degree less than $n$ with respect to the weight $x^{\beta+1}(1-x)^{\alpha-1}$, which proves (\ref{eq:lem1}). The more general expression
(\ref{eq:lemm}) follows by applying (\ref{eq:lem1}) $m$ times.
\qed

By using this lemma we arrive at
\begin{multline*}
   (n+m+\alpha_0+\alpha_1+1)_n   P_{n,m}^{(\alpha_0,\alpha_1,\alpha_2)}(x) \\
  = n!  \sum_{k=0}^n \binom{\alpha_1+n}{k} \binom{\alpha_0+m+n}{n-k}
   x^{n-k}(x-1)^{k} P_m^{(\alpha_0+k,\alpha_2+n-k)}(x). 
\end{multline*}
For the Jacobi polynomial we have the expansion
\begin{equation}  \label{eq:jacobiex}
  (\alpha+\beta+n+1)_n P_n^{(\alpha,\beta)}(x) =
 n! \sum_{j=0}^n \binom{\beta+n}{j} \binom{\alpha+n}{n-j}
 x^{n-j} (x-1)^{j},
\end{equation}
which can easily be obtained from the Rodrigues formula (\ref{eq:jacobiRod})
by using Leibniz' formula,
so that we finally find
\begin{multline}  \label{eq:pineirods}
   (n+m+\alpha_0+\alpha_1+1)_n (n+m+\alpha_0+\alpha_2+1)_m   P_{n,m}^{(\alpha_0,\alpha_1,\alpha_2)}(x) \\
  = n! m! \sum_{k=0}^n \sum_{j=0}^m \binom{\alpha_1+n}{k} \binom{\alpha_0+m+n}{n-k} \binom{\alpha_2+n+m-k}{j} 
\binom{\alpha_0+k+m}{m-j}
   x^{n+m-k-j}(x-1)^{k+j}. 
\end{multline}
We can explicitly find the first few coefficients of $P_{m,n}^{(\alpha_0,\alpha_1,\alpha_2)}(x)$ from this expression. We
introduce the notation
\begin{eqnarray*}
  K_{n,m} & = &  \frac{n! m!}{(n+m+\alpha_0+\alpha_1+1)_n (n+m+\alpha_0+\alpha_2+1)_m} \\
   & = & \binom{\alpha_0+\alpha_1+2n+m}{n}^{-1} \binom{\alpha_0+\alpha_2+2m+n}{m}^{-1}.
\end{eqnarray*}
First let us check that the polynomial is indeed monic by working out
the coefficient of $x^{m+n}$. This is given by
\[ K_{n,m} \sum_{k=0}^n \sum_{j=0}^m \binom{\alpha_1+n}{k} \binom{\alpha_0+m+n}{n-k} \binom{\alpha_2+n+m-k}{j} 
\binom{\alpha_0+k+m}{m-j}. \]
The sum over $j$ can be evaluated using the Chu-Vandermonde identity
\[ \sum_{j=0}^m  \binom{\alpha_2+n+m-k}{j} 
\binom{\alpha_0+k+m}{m-j} = \binom{\alpha_0+\alpha_2+n+2m}{m}, \]
which is independent of $k$.
The remaining sum over $k$ can also be evaluated and gives
\[  \sum_{k=0}^n \binom{\alpha_1+n}{k} \binom{\alpha_0+m+n}{n-k}
 = \binom{\alpha_0+\alpha_1+m+2n}{n}, \]
and the double sum is therefore equal to $K_{n,m}^{-1}$, showing
that this polynomial is indeed monic. Now let us write
\[  P_{n,m}^{(\alpha_0,\alpha_1,\alpha_2)}(x) =
  x^{m+n} + A_{n,m} x^{n+m-1} + B_{n,m}x^{n+m-2} + C_{n,m} x^{n+m-3} +
  \cdots . \]
The coefficient $A_{n,m}$ of $x^{m+n-1}$ is given by
\[ -K_{n,m} \sum_{k=0}^n \sum_{j=0}^m (k+j) \binom{\alpha_1+n}{k} \binom{\alpha_0+m+n}{n-k} \binom{\alpha_2+n+m-k}{j} 
\binom{\alpha_0+k+m}{m-j}. \]
This double sum can again be evaluated using  Chu-Vandermonde and
gives
\[  A_{n,m} = -\frac{n(\alpha_1+n)(\alpha_0+\alpha_2+n+m)+m(\alpha_2+n+m)
(\alpha_0+\alpha_1+2n+m)}{(\alpha_0+\alpha_1+2n+m)(\alpha_0+\alpha_2+n+2m)}.
\]
Similarly we can compute the coefficient $B_{n,m}$ of $x^{n+m-2}$ and
the coefficient
$C_{n,m}$ of $x^{m+n-3}$, but the computation is rather
lengthy. Once these coefficients have been determined, one can compute
the coefficients in the recurrence relation
\[  xP_{n}(x) = P_{n+1}(x) + b_n P_n(x) + c_n P_{n-1}(x) + d_n P_{n-2}(x), \]
where
\[  P_{2n}(x) = P_{n,n}^{(\alpha_0,\alpha_1,\alpha_2)}(x), \quad
   P_{2n+1}(x) = P_{n+1,n}^{(\alpha_0,\alpha_1,\alpha_2)}(x). \]
Indeed, by comparing coefficients we have
\begin{equation}  \label{eq:b}
  b_{2n} = A_{n,n}-A_{n+1,n}, \quad b_{2n+1} = A_{n+1,n}-A_{n+1,n+1},
\end{equation}
which gives
\begin{eqnarray*}
 b_{2n} & = & 
 [36{n}^{4}+ (48\alpha_0+28\alpha_1+20\alpha_2+38){n}^{3} \\
 && +\  (21\alpha_0^2+8\alpha_1^2+4\alpha_2^2+30\alpha_0\alpha_1+18\alpha_0\alpha_2
+15\alpha_1\alpha_2+39\alpha_0+19\alpha_1+19\alpha_2+9){n}^{2} \\
 && +\ (3\alpha_0^3+10\alpha_0^2\alpha_1+
  4\alpha_0^2\alpha_2 + 6\alpha_0\alpha_1^2+2\alpha_0\alpha_2^2+
11\alpha_0\alpha_1\alpha_2 + 5\alpha_1^2\alpha_2+3\alpha_1\alpha_2^2 \\
&& +\ 12\alpha_0^2+3\alpha_1^2+3\alpha_2^2+13\alpha_0\alpha_1+
  13\alpha_0\alpha_2+8\alpha_1\alpha_2+  
6\alpha_0+  3\alpha_1+3\alpha_2)n \\
 && +\ \alpha_0^2+\alpha_0\alpha_1+ \alpha_2\alpha_1^2+2\alpha_2\alpha_1^2\alpha_0+2\alpha_0^2\alpha_1+
\alpha_1^2\alpha_0+\alpha_2^2\alpha_0+\alpha_2^2\alpha_1+\alpha_0^3\alpha_1 \\ 
&& +\ \alpha_0^2\alpha_1^2+\alpha_2^2\alpha_0\alpha_1+\alpha_2^2\alpha_1^2
+2\alpha_2\alpha_0^2\alpha_1+3\alpha_2\alpha_1\alpha_0+2\alpha_2\alpha_0^2
+\alpha_1\alpha_2+\alpha_0^3+\alpha_0\alpha_2] \\
&& \times \ 
 (3n+\alpha_0+\alpha_2)^{-1} (3n+\alpha_0+\alpha_1)^{-1} (3n+\alpha_0+\alpha_2+1)^{-1}  (3n+\alpha_0+\alpha_1+2)^{-1} ,
\end{eqnarray*}
and
\begin{eqnarray*}
 b_{2n+1} & = & 
 [36{n}^{4}+ (48\alpha_0+20\alpha_1+28\alpha_2+106){n}^{3}\\
&& +\ (21\alpha_0^2+4\alpha_1^2+8\alpha_2^2+18\alpha_0\alpha_1
 +30\alpha_0\alpha_2+15\alpha_1\alpha_2+105\alpha_0+
41\alpha_1+65\alpha_2+111){n}^{2} \\
&& +\
 (3\alpha_0^3+4\alpha_0^2\alpha_1+10\alpha_0^2\alpha_2
  +2\alpha_0\alpha_1^2+6\alpha_0\alpha_2^2+11\alpha_0\alpha_1\alpha_2
 +3\alpha_1^2\alpha_2 + 5\alpha_1\alpha_2^2 \\
&& +\ 30\alpha_0^2 +5\alpha_1^2+13\alpha_2^2+23\alpha_0\alpha_1 +
 47\alpha_0\alpha_2+22\alpha_1\alpha_2+72\alpha_0+25\alpha_1+
  49\alpha_2+48)n \\
&& +\ 18\alpha_0\alpha_2+8\alpha_2\alpha_0^2+4\alpha_1+4\alpha_2^2\alpha_1
  +8\alpha_1\alpha_2+2\alpha_0^3+5\alpha_2^2\alpha_0+
   8\alpha_2\alpha_1\alpha_0+12\alpha_2 \\
&& +\ 7+15\alpha_0+\alpha_2^2\alpha_1^2+10\alpha_0^2+6\alpha_0\alpha_1
   +2\alpha_2\alpha_1^2+2\alpha_0^2\alpha_1+\alpha_1^2\alpha_0
   +5\alpha_2^2+\alpha_2\alpha_0^3 \\ 
&& +\ \alpha_2^2\alpha_0^2+\alpha_1^2+\alpha_2\alpha_1^2\alpha_0+
  2\alpha_2\alpha_0^2\alpha_1+2\alpha_2^2\alpha_0\alpha_1] \\
&& \times \
 (3n+\alpha_0+\alpha_2+1)^{-1} (3n+\alpha_0+\alpha_1+2 )^{-1}
 (3n+\alpha_0+\alpha_2+3)^{-1} (3n+\alpha_0+\alpha_1+3)^{-1} .
\end{eqnarray*}
For the recurrence coefficient $c_n$ we have the formulas
\begin{equation}  \label{eq:c}
 c_{2n} = B_{n,n} - B_{n+1,n} - b_{2n} A_{n,n}, \quad
    c_{2n+1} = B_{n+1,n} - B_{n+1,n+1} - b_{2n+1} A_{n+1,n}, 
\end{equation}
which, after some computation (and using Maple V), gives
\begin{eqnarray*}
c_{2n}& = & 
 n(2n+\alpha_0)(2n+\alpha_0+\alpha_1)(2n+\alpha_0+\alpha_2)
  \\
&& \times\ [ 54n^4 + (63\alpha_0 + 45\alpha_1 + 45 \alpha_2)n^3 \\
&& +\ (24\alpha_0^2 + 8 \alpha_1^2+ 8\alpha_2^2 + 42 \alpha_0\alpha_1
  + 42 \alpha_0\alpha_2 + 44 \alpha_1\alpha_2 - 8) n^2 \\
&& +\ (3 \alpha_0^3 + \alpha_1^3 + \alpha_2^3 + 12 \alpha_0^2\alpha_1
  + 12 \alpha_0^2\alpha_2 + 3 \alpha_0\alpha_1^2 + 3 \alpha_0\alpha_2^2
 + 33 \alpha_0\alpha_1\alpha_2 +  8\alpha_1^2\alpha_2 \\
&& +\ 8\alpha_1\alpha_2^2 -3\alpha_0 - 4 \alpha_1
  -4\alpha_2) n \\
 && +\ \alpha_0^3\alpha_1+\alpha_0^3\alpha_2 + 6 \alpha_0^2\alpha_1\alpha_2
  + \alpha_1^3\alpha_2 + \alpha_1\alpha_2^3 
+ 3 \alpha_0 \alpha_1^2\alpha_2
  + 3 \alpha_0\alpha_1\alpha_2^2 - \alpha_0\alpha_1 - \alpha_0\alpha_2
 - 2 \alpha_1\alpha_2 ]  \\
&& \times\ (3n+\alpha_0+\alpha_1+1)^{-1}
(3n+\alpha_0+\alpha_2+1)^{-1} (3n+\alpha_0+\alpha_1)^{-2}
  (3n+\alpha_0+\alpha_2)^{-2} \\
&& (3n+\alpha_0+\alpha_1-1)^{-1}(3n+\alpha_0+\alpha_2-1)^{-1}
\end{eqnarray*}
and
\begin{eqnarray*}
c_{2n+1}& = & 
 (2n+\alpha_0+1)(2n+\alpha_0+\alpha_1+1)(2n+\alpha_0+\alpha_2+1)
  \\
&& \times\ [ 54n^5 + (63\alpha_0 + 45\alpha_1 + 45 \alpha_2+135)n^4 \\
&& +\ (24\alpha_0^2 + 8 \alpha_1^2+ 8\alpha_2^2 + 42 \alpha_0\alpha_1
  + 42 \alpha_0\alpha_2 + 44 \alpha_1\alpha_2 - +126\alpha_0 +76\alpha_1
 +104\alpha_2 + 120) n^3 \\
&& +\ (3 \alpha_0^3 + \alpha_1^3 + \alpha_2^3 + 12 \alpha_0^2\alpha_1
  + 12 \alpha_0^2\alpha_2 + 3 \alpha_0\alpha_1^2 + 3 \alpha_0\alpha_2^2
 + 33 \alpha_0\alpha_1\alpha_2 + 8\alpha_1^2\alpha_2 \\
&& +\ 8\alpha_1\alpha_2^2 + 36\alpha_0^2+5\alpha_1^2
+19\alpha_2^2 + +54\alpha_0\alpha_1 + 72\alpha_0\alpha_2 + 66 \alpha_1\alpha_2 + 87 \alpha_0 + 39 \alpha_1 \\
&& +\ 81 \alpha_2 + 45) n^2 \\ 
&& +\ ( \alpha_0^3\alpha_1+\alpha_0^3\alpha_2 + 6 \alpha_0^2\alpha_1\alpha_2
  + \alpha_1^3\alpha_2 + \alpha_1\alpha_2^3  + 3 \alpha_0 \alpha_1^2\alpha_2
 + 3 \alpha_0\alpha_1\alpha_2^2  + 3\alpha_0^3 + 2\alpha_2^3 \\ 
&& +\ 12\alpha_0^2\alpha_1 + 12\alpha_0^2\alpha_2 + 6\alpha_0\alpha_2^2
 + 33 \alpha_0\alpha_1\alpha_2 
+ 5\alpha_1^2\alpha_2 +11 \alpha_1\alpha_2^2
 +18 \alpha_0^2 +20\alpha_0\alpha_1 \\ 
&& +\ 38 \alpha_0\alpha_2 + 14\alpha_2^2
 + 26\alpha_1\alpha_2 + 24 \alpha_0 + 6\alpha_1 +24\alpha_2 +6)n \\
&& +\ \alpha_0^3\alpha_1 + 3\alpha_0^2\alpha_1\alpha_2 + 3\alpha_0\alpha_1\alpha_2^2 + \alpha_1\alpha_2^3 +\alpha_0^3 + \alpha_2^3
 +3\alpha_0^2\alpha_1 + 3\alpha_0^2\alpha_2 + 6\alpha_0\alpha_1\alpha_2 \\
&& +\ 3\alpha_0\alpha_2^2 + 3\alpha_1\alpha_2^2 + 3\alpha_0^2 + 3\alpha_2^2 
 + 2\alpha_0\alpha_1 + 6\alpha_0\alpha_2 + 2\alpha_1\alpha_2 +2\alpha_0
+ 2\alpha_2  ]  \\
&& \times\ (3n+\alpha_0+\alpha_1+3)^{-1}
(3n+\alpha_0+\alpha_2+2)^{-1} (3n+\alpha_0+\alpha_1+2)^{-2}
  (3n+\alpha_0+\alpha_2+1)^{-2} \\
&& (3n+\alpha_0+\alpha_1+1)^{-1}(3n+\alpha_0+\alpha_2)^{-1}.
\end{eqnarray*}
Finally, for $d_n$ we have 
\begin{eqnarray}  \label{eq:d}
 d_{2n} & = & C_{n,n} - C_{n+1,n} - b_{2n} B_{n,n} - c_{2n} A_{n,n-1} ,
  \nonumber \\
  d_{2n+1} & = & C_{n+1,n} - C_{n+1,n+1} - b_{2n+1} B_{n+1,n} - c_{2n+1} A_{n,n},
\end{eqnarray}
giving
\begin{eqnarray*}
  d_{2n} & = & n(2n+\alpha_0)(2n+\alpha_0-1)(2n+\alpha_0+\alpha_1)
   (2n+\alpha_0+\alpha_1-1) \\
 &&  (2n+\alpha_0+\alpha_2)(2n+\alpha_0+\alpha_2-1)
  (n+\alpha_1)(n+\alpha_1-\alpha_2) \\
 && (3n+1+\alpha_0+\alpha_1)^{-1}(3n+\alpha_0+\alpha_1)^{-2}
    (3n+\alpha_0+\alpha_2)^{-1}(3n-1+\alpha_0+\alpha_1)^{-2} \\
 &&   (3n-1+\alpha_0+\alpha_2)^{-1} (3n-2+\alpha_0+\alpha_1)^{-1}
    (3n-2+\alpha_0+\alpha_2)^{-1}
\end{eqnarray*}
and
\begin{eqnarray*}
  d_{2n+1} & = & n(2n+1+\alpha_0)(2n+\alpha_0)(2n+\alpha_0+\alpha_1)
   (2n+1+\alpha_0+\alpha_1) \\
 &&  (2n+1+\alpha_0+\alpha_2)(2n+\alpha_0+\alpha_2)
  (n+\alpha_2)(n+\alpha_2-\alpha_1) \\
 && (3n+2+\alpha_0+\alpha_1)^{-1}(3n+2+\alpha_0+\alpha_2)^{-1}
    (3n+1+\alpha_0+\alpha_1)^{-1}(3n+1+\alpha_0+\alpha_2)^{-2} \\
 &&   (3n+\alpha_0+\alpha_1)^{-1} (3n+\alpha_0+\alpha_2)^{-2}
    (3n-1+\alpha_0+\alpha_2)^{-1}.
\end{eqnarray*}
These formulas are rather lengthy, but explicit knowledge of them will be
useful in what follows. Observe that for large $n$ we have
\begin{eqnarray*}
  \lim_{n \to \infty} b_n & = & \frac49 = 3 \left( \frac4{27} \right) , \\
  \lim_{n \to \infty} c_n & = & \frac{16}{243} = 3 \left( \frac4{27} \right)^2 , \\
  \lim_{n \to \infty} d_n & = & \frac{64}{19683} = \left( \frac{4}{27} \right)^3.
\end{eqnarray*}

\subsection{Multiple Laguerre polynomials (first kind)}
In the same spirit as for the Jacobi-Pi\~neiro polynomials, we can
consider two different families of multiple Laguerre polynomials.
The \textit{multiple Laguerre polynomials of the first kind}
$L_{\vec{n}}^{\vec{\alpha}}(x)$ are orthogonal on $[0,\infty)$
with respect to the $r$ weights $w_j(x) = x^{\alpha_j} e^{-x}$,
where $\alpha_j > -1$ for $j=1,2,\ldots,r$. So these weights have the
same exponential decrease at $\infty$ but have different singularities
at $0$. Again we assume $\alpha_i -\alpha_j \notin \mathbb{Z}$ in order
to have an AT system. 
These polynomials were first considered by Sorokin \cite{sor2},
\cite{sor6}.
The raising operators are given by
\begin{equation}
 \frac{d}{dx} \left( x^{\alpha_j} e^{-x} L_{\vec{n}}^{\vec{\alpha}}(x)
 \right) = - x^{\alpha_j-1} e^{-x} 
   L_{\vec{n}+\vec{e}_j}^{\vec{\alpha}-\vec{e}_j}(x), \qquad
  j=1,\ldots,r, 
\end{equation}
and a repeated application of these operators gives the Rodrigues formula
\begin{equation}  \label{eq:lagIRod}
 (-1)^{|\vec{n}|}   L_{\vec{n}}^{\vec{\alpha}}(x)
 = e^x \prod_{j=1}^r \left[ x^{-\alpha_j} 
  \frac{d^{n_j}}{dx^{n_j}} x^{n_j+\alpha_j} \right] e^{-x}.
\end{equation}

When $r=2$ one can use this Rodrigues formula to obtain an explicit expression for these multiple Laguerre polynomials, from which one can
compute the recurrence coefficients in
\[ xP_n(x) = P_{n+1}(x) + b_n P_n(x) + c_n P_{n-1}(x) + d_n P_{n-2}(x), \]
where $P_{2n}(x) = L_{n,n}^{(\alpha_1,\alpha_2)}(x)$ and
$P_{2n+1}(x) = L_{n+1,n}^{(\alpha_1,\alpha_2)}(x)$. But having done
all that work for Jacobi-Pi\~neiro polynomials, it is much easier to use
the limit relation 
\begin{equation}
  L_{n,m}^{(\alpha_1,\alpha_2)}(x) = \lim_{\alpha_0 \to \infty}
   \alpha_0^{n+m} P_{n,m}^{(\alpha_0,\alpha_1,\alpha_2)}(x/\alpha_0).
\end{equation}
The recurrence coefficients can then be found in terms of the following
limits of the corresponding recurrence coefficients of Jacobi-Pi\~neiro
polynomials:
\begin{eqnarray*}
    b_n & =&  \lim_{\alpha_0 \to \infty} b_n^{(\alpha_0,\alpha_1,\alpha_2)}     \alpha_0,  \\
    c_n & =&  \lim_{\alpha_0 \to \infty} c_n^{(\alpha_0,\alpha_1,\alpha_2)}     \alpha_0^2 , \\
   d_n & =&  \lim_{\alpha_0 \to \infty} d_n^{(\alpha_0,\alpha_1,\alpha_2)}     \alpha_0^3 ,
\end{eqnarray*}
giving
\begin{eqnarray*}
 b_{2n} & = & 3n+\alpha_1+1, \\
 b_{2n+1} & = & 3n+\alpha_2+2, \\
 c_{2n} & = & n(3n +\alpha_1+\alpha_2), \\
 c_{2n+1} & = & 3n^2+(\alpha_1+\alpha_2+3)n + \alpha_1+1, \\
 d_{2n} & = & n(n+\alpha_1)(n+\alpha_1-\alpha_2), \\
 d_{2n+1} & = & n (n+\alpha_2)(n+\alpha_2-\alpha_1).
\end{eqnarray*}
Observe that for large $n$ we have
\begin{eqnarray*}
	\lim_{n \to \infty} \frac{b_n}n & = & \frac32 = 3 \left( \frac12 \right), \\
      \lim_{n \to \infty} \frac{c_n}{n^2} & = & \frac34 = 3 \left( \frac12 \right)^2 , \\
   \lim_{n \to \infty} \frac{d_n}{n^3} & = & \frac18 =  \left( \frac12 \right)^3.
\end{eqnarray*}

\subsection{Multiple Laguerre polynomials (second kind)}
Another family of multiple Laguerre polynomials is given by the weights
$w_j(x) = x^{\alpha_0} e^{-c_j x}$ on $[0,\infty)$, with $c_j > 0$
and $c_i\neq c_j$ for $i\neq j$. So now the weights have the same
singularity at the origin but different exponential rates at infinity.
These \textit{multiple Laguerre polynomials of the second kind} 
$L_{\vec{n}}^{(\alpha_0,\vec{c})}(x)$ appear
already in \cite[p.~160]{nikisor}. The raising operators are
\begin{equation}
 \frac{d}{dx} \left( x^{\alpha_0} e^{-c_jx} L_{\vec{n}}^{(\alpha_0,\vec{c})}(x)
 \right) = - c_j x^{\alpha_0-1} e^{-c_jx} 
   L_{\vec{n}+\vec{e}_j}^{(\alpha_0-1,\vec{c})}(x), \qquad
  j=1,\ldots,r, 
\end{equation}
and a repeated application of these operators gives the Rodrigues formula
\begin{equation}  \label{eq:lagIIRod}
 (-1)^{|\vec{n}|}   \prod_{j=1}^r c_j^{n_j} \
L_{\vec{n}}^{(\alpha_0,\vec{c})}(x)
 = x^{-\alpha_0} \prod_{j=1}^r \left[ e^{c_j x} 
  \frac{d^{n_j}}{dx^{n_j}} e^{-c_j x} \right] x^{|\vec{n}|+\alpha_0}.
\end{equation}
These polynomials are also a limit case of the Jacobi-Pi\~neiro polynomials.
For the case $r=2$ we have
\begin{equation}
  L_{n,m}^{(\alpha_0,c_1,c_2)}(x) = \lim_{\alpha \to \infty}
   (-\alpha)^{n+m} P_{n,m}^{(\alpha_0,c_1\alpha,c_2\alpha)}(1-x/\alpha).
\end{equation}
The recurrence coefficients can be obtained from the corresponding recurrence coefficients of Jacobi-Pi\~neiro polynomials by
\begin{eqnarray*}
    b_n & =&  \lim_{\alpha \to \infty} 
   (1-b_n^{(\alpha_0,c_1\alpha,c_2\alpha)}) \alpha,  \\
    c_n & =&  \lim_{\alpha \to \infty} 
   c_n^{(\alpha_0,c_1\alpha,c_2\alpha)} \alpha^2 , \\
   d_n & =&  \lim_{\alpha \to \infty} 
   -d_n^{(\alpha_0,c_1\alpha,c_2\alpha)} \alpha^3 ,
\end{eqnarray*}
giving
\begin{eqnarray*}
  b_{2n} & = & \frac{n(c_1+3c_2)+c_2+\alpha_0c_2}{c_1c_2}, \\
  b_{2n+1} & = & \frac{n(3c_1+c_2)+2c_1+c_2+\alpha_0c_1}{c_1c_2}, \\
  c_{2n} & = & \frac{n(2n+\alpha_0)(c_1^2+c_2^2)}{c_1^2c_2^2}, \\
  c_{2n+1} & = & \frac{2n^2(c_1^2+c_2^2) + n[c_1^2+3c_2^2+\alpha_0(c_1^2+c_2^2)]+c_2^2+\alpha_0c_2^2}{c_1^2c_2^2}, \\
  d_{2n} & = & \frac{n(2n+\alpha_0)(2n+\alpha_0-1)(c_2-c_1)}{c_1^3c_2}, \\
  d_{2n+1} & = & \frac{n(2n+\alpha_0)(2n+\alpha_0+1)(c_1-c_2)}{c_1c_2^3}.
\end{eqnarray*}
Observe that for large $n$ we have
\begin{eqnarray*}
	\lim_{n \to \infty} \frac{b_n}n & = &
   \begin{cases} \displaystyle
       \frac{c_1+3c_2}{2c_1c_2} & \mathrm{if\ } n \equiv 0 \pmod{2} , \\
               \displaystyle
       \frac{3c_1+c_2}{2c_1c_2} & \mathrm{if\ } n \equiv 1 \pmod{2} ,
   \end{cases} \\
      \lim_{n \to \infty} \frac{c_n}{n^2} & = & 
   \frac{c_1^2+c_2^2}{2c_1^2c_2^2}, \\
   \lim_{n \to \infty} \frac{d_n}{n^3} & = &
     \begin{cases} \displaystyle
       \frac{c_2-c_1}{2c_1^3c_2} & \mathrm{if\ } n \equiv 0 \pmod{2} , \\
          \displaystyle
       \frac{c_1-c_2}{2c_1c_2^3} & \mathrm{if\ } n \equiv 1 \pmod{2} .
   \end{cases}
\end{eqnarray*}

\subsection{Multiple Hermite polynomials}
Finally we can consider the weights $w_j(x) = e^{-x^2+c_j x}$ on $(-\infty,\infty)$, for $j=1,2,\ldots,r$ and $c_j$ different real numbers.
The \textit{multiple
Hermite polynomials} $H_{\vec{n}}^{\vec{c}}(x)$ once more have
raising operators and a Rodrigues formula, and they are also limiting cases
of the Jacobi-Pi\~neiro polynomials, but also of the multiple
Laguerre polynomials of the second kind. For $r=2$ this is
\begin{equation}
  H_{n,m}^{(c_1,c_2)}(x) = \lim_{\alpha \to \infty}
   (2\sqrt{\alpha})^{n+m} P_{n,m}^{(\alpha,\alpha + c_1\sqrt{\alpha},
 \alpha+c_2\sqrt{\alpha})}\left( \frac{x+\sqrt{\alpha}}{2\sqrt{\alpha}}
  \right) ,
\end{equation}
so that the recurrence coefficients can be obtained from the
Jacobi-Pi\~neiro case by
\begin{eqnarray*}
    b_n & =&  \lim_{\alpha \to \infty} 
   2(b_n^{(\alpha,\alpha+c_1\sqrt{\alpha},\alpha+c_2\sqrt{\alpha})}-\frac12) \sqrt{\alpha},  \\
    c_n & =&  \lim_{\alpha \to \infty} 
   4 c_n^{(\alpha,\alpha+c_1\sqrt{\alpha},\alpha+c_2\sqrt{\alpha})} \alpha , \\
   d_n & =&  \lim_{\alpha \to \infty} 
   8 d_n^{(\alpha,\alpha+c_1\sqrt{\alpha},\alpha+c_2\sqrt{\alpha})} 
  (\sqrt{\alpha})^3 .
\end{eqnarray*}
This gives
\begin{eqnarray*}
  b_{2n} & = & c_1/2 , \\
  b_{2n+1} & = & c_2/2, \\
  c_n  & = & n/2 , \\
  d_{2n} & = & n(c_1-c_2)/4 , \\
  d_{2n+1} & = & n (c_2-c_1)/4. 
\end{eqnarray*}
Alternatively, we can use the limit transition from the multiple
Laguerre polynomials of the first kind: 
\begin{equation}
  H_{n,m}^{(c_1,c_2)}(x) = \lim_{\alpha \to \infty}
   \alpha^{n+m} L_{n,m}^{(\alpha + c_1\sqrt{\alpha/2},
 \alpha+c_2\sqrt{\alpha/2})}(\sqrt{2\alpha}x+\alpha) .
\end{equation}
The  recurrence coefficients are then also given in terms of the following
limits of the recurrence coefficients of the multiple Laguerre polynomials
of the first kind 
\begin{eqnarray*}
    b_n & =&  \lim_{\alpha \to \infty} 
   (b_n^{(\alpha+c_1\sqrt{\alpha/2},\alpha+c_2\sqrt{\alpha/2})}-\alpha)/ \sqrt{2\alpha},  \\
    c_n & =&  \lim_{\alpha \to \infty} 
    c_n^{(\alpha+c_1\sqrt{\alpha/2},\alpha+c_2\sqrt{\alpha/2})}/(2\alpha) , \\
   d_n & =&  \lim_{\alpha \to \infty} 
    d_n^{(\alpha+c_1\sqrt{\alpha/2},\alpha+c_2\sqrt{\alpha/2})} / 
  (\sqrt{2\alpha})^3 ,
\end{eqnarray*}
which leads to the same result.
Observe that for large $n$ we have
\begin{eqnarray*}
  \lim_{n \to \infty}  \frac{b_n}{\sqrt{n}} & = & 0 , \\
   \lim_{n \to \infty}  \frac{c_n}{n} & = & \frac12 , \\
  \lim_{n \to \infty}  \frac{d_n}{(\sqrt{n})^3} & = & 0 .
\end{eqnarray*}

\subsection{Jacobi-Angelesco polynomials}
The following system is probably the first that was investigated in detail (\cite{kalya1}, \cite{kalyaron}). 
It is an Angelesco system with weights
$w_1(x) = |h(x)|$ on $[a,0]$ (with $a < 0$) and $w_2(x)= |h(x)|$ on $[0,1]$, where
$h(x) = (x-a)^\alpha x^\beta (1-x)^\gamma$ and $\alpha, \beta , \gamma > -1$. Hence the same weight is used
for both weights $w_1$ and $w_2$ but on two touching intervals.  The Jacobi-Angelesco polynomials $P_{n,m}^{(\alpha,\beta,\gamma)}(x;a)$ therefore
satisfy the orthogonality relations
\begin{eqnarray}
  \int_a^0 P_{n,m}^{(\alpha,\beta,\gamma)}(x;a) (x-a)^\alpha |x|^\beta
   (1-x)^\gamma x^k\, dx & = & 0, \qquad k = 0,1,2,\ldots,n-1, \\
  \int_0^1 P_{n,m}^{(\alpha,\beta,\gamma)}(x;a) (x-a)^\alpha x^\beta
   (1-x)^\gamma x^k\, dx & = & 0, \qquad k = 0,1,2,\ldots,m-1.
\end{eqnarray}
The function $h(x)$ satisfies a Pearson equation
\[   (x-a)x(1-x) h'(x) = [\alpha x(1-x) +\beta (x-a)(1-x) - \gamma (x-a)x]
    h(x), \]
where $(x-a)x(1-x)$ is now a polynomial of degree 3. Using this relation,
we can integrate the orthogonality relations by part to see that
\begin{multline}   
  \frac{d}{dx} \left[ (x-a)^\alpha x^\beta (1-x)^\gamma  
    P_{n,m}^{(\alpha,\beta,\gamma)}(x;a) \right] \\
   = -(\alpha+\beta+\gamma+n+m) (x-a)^{\alpha-1} x^{\beta-1}
   (1-x)^{\gamma-1} P_{n+1,m+1}^{(\alpha-1,\beta-1,\gamma-1)}(x;a), 
\end{multline}
which raises both indices of the multi-index $(n,m)$.
Repeated use of this raising operation gives the Rodrigues formula
\begin{multline}  \label{eq:jacangRod}
  \frac{d^m}{dx^m} \left[ (x-a)^{\alpha+m} x^{\beta+m} (1-x)^{\gamma+m}
    P_{k,0}^{(\alpha+m,\beta+m,\gamma+m)}(x;a) \right] \\
  = (-1)^m (\alpha+\beta+\gamma+k+2m+1)_m (x-a)^\alpha x^\beta (1-x)^\gamma
   P_{m+k,m}^{(\alpha,\beta,\gamma)}(x;a). 
\end{multline}
For $k=0$ and $m=n$, this then gives 
\begin{multline}  \label{eq:jacangRodn}
  \frac{d^n}{dx^n} \left[ (x-a)^{\alpha+n} x^{\beta+n} (1-x)^{\gamma+n}
     \right] \\
  = (-1)^n (\alpha+\beta+\gamma+2n+1)_n (x-a)^\alpha x^\beta (1-x)^\gamma
   P_{n,n}^{(\alpha,\beta,\gamma)}(x;a). 
\end{multline}
Use Leibniz' formula to find
\begin{multline*}
 (-1)^n (\alpha+\beta+\gamma+2n+1)_n (x-a)^\alpha x^\beta (1-x)^\gamma
   P_{n,n}^{(\alpha,\beta,\gamma)}(x;a) \\
= \sum_{k=0}^n \binom{n}{k}
  \left( \frac{d^k}{dx^k} x^{\beta+n}(1-x)^{\gamma+n} \right) \ 
 \left( \frac{d^{n-k}}{dx^{n-k}} (x-a)^{\alpha+n} \right).
\end{multline*}
Now use the Rodrigues formula for the Jacobi polynomials
(\ref{eq:jacobiRod}) to find
\begin{multline*}
  \binom{\alpha+\beta+\gamma+3n}{n} P_{n,n}^{(\alpha,\beta,\gamma)}(x;a)  \\  =  \sum_{k=0}^n (-1)^{n-k}  \binom{\beta+\gamma+2n}{k}\binom{\alpha+n}{n-k}
  (x-a)^kx^{n-k}(1-x)^{n-k} P_{k}^{(\gamma+n-k,\beta+n-k)}(x).
\end{multline*}
Use of the expansion (\ref{eq:jacobiex}) for the Jacobi polynomial gives
\begin{eqnarray} 
   \lefteqn{\binom{\alpha+\beta+\gamma+3n}{n} P_{n,n}^{(\alpha,\beta,\gamma)}(x;a)} & &  \nonumber \\
& = & \sum_{k=0}^n\sum_{j=0}^k 
   \binom{\alpha+n}{n-k}
   \binom{\beta+n}{j} \binom{\gamma+n}{k-j}
  (x-a)^k x^{n-j}(x-1)^{n-k+j}   \label{eq:jacang1}  \\
 & = & \sum_{k=0}^n\sum_{j=0}^{n-k} 
   \binom{\alpha+n}{k}
   \binom{\beta+n}{j} \binom{\gamma+n}{n-k-j}
  (x-a)^{n-k} x^{n-j}(x-1)^{k+j},  \label{eq:jacang2}
\end{eqnarray}
where the last equation follows by the change of variable $k \mapsto n-k$.
If we write this in terms of Pochhammer symbols, then
\begin{eqnarray*} 
   \lefteqn{\binom{\alpha+\beta+\gamma+3n}{n} P_{n,n}^{(\alpha,\beta,\gamma)}(x;a)} & &  \nonumber \\
& = &  \frac{(\gamma+1)_n}{n!}
 \sum_{k=0}^n\sum_{j=0}^{n-k}  \frac{(-n)_{k+j}(-\alpha-n)_k(-\beta-n)_j}
    {(\gamma+1)_{k+j}k!j!} (x-a)^{n-k} (x-1)^{k+j} x^{n-j} \\
 & = & x^n(x-a)^n \binom{\gamma+n}{n} F_1\left(-n,-\alpha-n,-\beta-n,\gamma+1; \frac{x-1}{x-a}, \frac{x-1}{x} \right),
\end{eqnarray*}
where
\[  F_1(a,b,b',c;x,y) = \sum_{m=0}^\infty \sum_{n=0}^\infty
   \frac{(a)_{m+n} (b)_m (b')_n}{(c)_{m+n}} \frac{x^m y^n}{m! n!} \]
is the first of Appell's hypergeometric functions of two variables.

For the polynomial $P_{n+1,n}^{(\alpha,\beta,\gamma)}(x;a)$ we have
the Rodrigues formula
\begin{multline}  \label{eq:jacangRodnp}
  \frac{d^n}{dx^n} \left[ (x-a)^{\alpha+n} x^{\beta+n} (1-x)^{\gamma+n}
     P_{1,0}^{(\alpha+n,\beta+n,\gamma+n)}(x;a) \right] \\
  = (-1)^n (\alpha+\beta+\gamma+2n+2)_n (x-a)^\alpha x^\beta (1-x)^\gamma
   P_{n+1,n}^{(\alpha,\beta,\gamma)}(x;a), 
\end{multline}
where $P_{1,0}^{(\alpha+n,\beta+n,\gamma+n)}(x;a) = x-X_n^{(\alpha,\beta,\gamma)}$
is the monic orthogonal polynomial of first degree for the
weight $(x-a)^{\alpha+n} |x|^{\beta+n} (1-x)^{\gamma+n}$ on $[a,0]$. If we write
down the orthogonality of this polynomial to the constant function,
\[  \int_a^0 (x-X_n^{(\alpha,\beta,\gamma)}) 
 (x-a)^{\alpha+n} |x|^{\beta+n} (1-x)^{\gamma+n} \, dx = 0, \]
then we see that
\[  X_n^{(\alpha,\beta,\gamma)} = 
   \frac{ \int_a^0 x(x-a)^{\alpha+n} |x|^{\beta+n} (1-x)^{\gamma+n} \, dx}
   {\int_a^0 (x-a)^{\alpha+n} |x|^{\beta+n} (1-x)^{\gamma+n} \, dx} . \]
A standard saddle point method gives the asymptotic behavior
\begin{equation}  \label{eq:Xinf}
  \lim_{n \to \infty} X_n^{(\alpha,\beta,\gamma)} = x_1,
\end{equation}
where $x_1$ is the zero of $g'(x)$ in $[a,0]$, where $g(x)=(x-a)x(1-x)$.
. Combining the Rodrigues
equation in (\ref{eq:jacangRodnp}) with the Rodrigues equation
(\ref{eq:jacangRodn}) shows that
\begin{equation} \label{eq:jacangnp}
   P_{n+1,n}^{(\alpha,\beta,\gamma)}(x;a)
     = x P_{n,n}^{(\alpha,\beta+1,\gamma)}(x;a) - X_n^{(\alpha,\beta,\gamma)} \frac{\alpha+\beta+\gamma+2n+1}{\alpha+\beta+\gamma+3n+1}
  P_{n,n}^{(\alpha,\beta,\gamma)}(x;a). 
\end{equation}

In order to compute the coefficients of the recurrence relation
\[  xP_n(x) = P_{n+1}(x) + b_n P_n(x) + c_n P_{n-1}(x) + d_n P_{n-2}(x), \]
where
\[  P_{2n}(x) = P_{n,n}^{(\alpha,\beta,\gamma)}(x;a), \quad
   P_{2n+1}(x) = P_{n+1,n}^{(\alpha,\beta,\gamma)}(x;a), \]
we will compute the first few coefficients of the polynomials
\[   P_{n,m}^{(\alpha,\beta,\gamma)}(x;a) = x^{m+n} +
    A_{n,m}x^{n+m-1} + B_{n,m} x^{m+n-2} + C_{n,m} x^{n+m-3} + \cdots . \]
First we take $n=m$. In order to check that our polynomial is monic, we
see from (\ref{eq:jacang1}) that the leading coefficient is given by
\[  \binom{\alpha+\beta+\gamma+3n}{n}^{-1}
  \sum_{k=0}^n \sum_{j=0}^k \binom{\alpha+n}{n-k}
   \binom{\beta+n}{j} \binom{\gamma+n}{k-j}. \]
Chu-Vandermonde gives
\[  \sum_{j=0}^k \binom{\beta+n}{j} \binom{\gamma+n}{k-j}
    = \binom{\beta+\gamma+2n}{k}, \]
and also
\[ \sum_{k=0}^n \binom{\alpha+n}{n-k}\binom{\beta+\gamma+2n}{k}
   = \binom{\alpha+\beta+\gamma+3n}{n}, \]
so that the leading coefficient is indeed 1. The coefficient $A_{n,n}$
of $x^{2n-1}$ is equal to
\[ - \binom{\alpha+\beta+\gamma+3n}{n}^{-1} \sum_{k=0}^n \sum_{j=0}^k
    \binom{\alpha+n}{n-k} \binom{\beta+n}{j} \binom{\gamma+n}{k-j}
 (ak+n-k+j) . \]
Working out this double sum gives
\begin{equation}  \label{eq:Ann}
  A_{n,n}^{(\alpha,\beta,\gamma)} =   \frac{-n[\alpha+\beta+2n +
  a(\beta+\gamma+2n)]}{\alpha+\beta+\gamma+3n}.
\end{equation}
For $P_{n+1,n}^{(\alpha,\beta,\gamma)}(x;a)$ the coefficient $A_{n+1,n}$
of $x^{2n}$ can be obtained from (\ref{eq:jacangnp})
\begin{equation}
  A_{n+1,n}^{(\alpha,\beta,\gamma)} = A_{n,n}^{(\alpha,\beta+1,\gamma)}
  - X_n^{(\alpha,\beta,\gamma)}
   \frac{\alpha+\beta+\gamma+2n+1}{\alpha+\beta+\gamma+3n+1}. 
\end{equation}
The coefficient $b_n$ in the recurrence relation can now be found from
(\ref{eq:b})
\begin{eqnarray*}
   b_{2n} & = & \frac{n[n+\gamma + a(n+\alpha)]}{(\alpha+\beta+\gamma
     +3n)(\alpha+\beta+\gamma+3n+1)} 
    + X_n^{(\alpha,\beta,\gamma)} 
   \frac{2n+\alpha+\beta+\gamma+1}{3n+\alpha+\beta+\gamma+1}, \\
  b_{2n+1} & = & \left( \rule{0pt}{2.5ex} 5n^2 + (4\alpha+4\beta+3\gamma+7)n +(\alpha+\beta+\gamma+1)(\alpha+\beta+2)  \right. \\
  & & +\ \left. \rule{0pt}{2.5ex} a[ 5n^2 + (3\alpha+4\beta+4\gamma+7)n
  + (\alpha+\beta+\gamma+1)(\beta+\gamma+2)]
   \right) \\
&& \times\  (\alpha+\beta+\gamma+3n+1)^{-1}(\alpha+\beta+\gamma+3n+3)^{-1} \\
  & & - X_n^{(\alpha,\beta,\gamma)} 
   \frac{2n+\alpha+\beta+\gamma+1}{3n+\alpha+\beta+\gamma+1}.
\end{eqnarray*}
The coefficient $B_{n,n}$ of $x^{2n-2}$ in
$P_{n,n}^{(\alpha,\beta,\gamma)}(x;a)$ is given by
\begin{eqnarray*}
    B_{n,n}^{(\alpha,\beta,\gamma)} & = & \frac{an(\alpha+\beta+\gamma+2n)(\beta+n)}{(\alpha+\beta+\gamma+3n)
 (\alpha+\beta+\gamma+3n-1)} \\
 && +\   \frac{n(n-1)}{2(\alpha+\beta+\gamma+3n)(\alpha+\beta+\gamma+3n-1)} \\
 && \times\ \left[ (\alpha+\beta+2n)(\alpha+\beta+2n-1) +
   2a(\alpha+\beta+2n)(\beta+\gamma+2n) \right. \\
  && \ \left.  +\ a^2 (\beta+\gamma+2n)(
  (\beta+\gamma+2n-1) \right],
\end{eqnarray*}
and from (\ref{eq:jacangnp}) we also find
\[  B_{n+1,n}^{(\alpha,\beta,\gamma)} = B_{n,n}^{(\alpha,\beta+1,\gamma)}
  - X_n^{(\alpha,\beta,\gamma)} A_{n,n}^{(\alpha,\beta,\gamma)}
   \frac{\alpha+\beta+\gamma+2n+1}{\alpha+\beta+\gamma+3n+1} . \]
Using  (\ref{eq:c}) then gives
\begin{eqnarray*}
 c_{2n} & = & \frac{n(\alpha+\beta+\gamma+2n)}{(\alpha+\beta+\gamma+3n-1)(\alpha+\beta+\gamma+3n)^2
    (\alpha+\beta+\gamma+3n-1)} \\
  & &  \times\ \left( \rule{0pt}{2.5ex}
 (\alpha+\beta+2n)(\gamma+n) - 2a (\alpha+n)(\gamma+n)  
 + a^2(\beta+\gamma+2n)(\alpha+n) \right),
\end{eqnarray*}
and  
\begin{eqnarray*}
  c_{2n+1} & = &    \frac{\alpha+\beta+\gamma+2n+1}{(\alpha+\beta+\gamma+3n+3)
  (\alpha+\beta+\gamma+3n+2)(\alpha+\beta+\gamma+3n+1)^2
 (\alpha+\beta+\gamma+3n)} \\
 & & \ \times\ \left( \rule{0pt}{2.5ex}
   n(n+\gamma)(\alpha+\beta+2n+1)(\alpha+\beta+\gamma+3n+3) \right. \\
 & & \ \ -\ a[24n^4+(29\alpha+41\beta+29\gamma+48)n^3    \\
 & & \ \ +\ (10\alpha^2+39\alpha\beta + 26\alpha\gamma +29\beta^2
    +39\beta\gamma+10\gamma^2 + 44\alpha +62\beta +44\gamma+30)n^2 \\
 & & \ \ +\ (\alpha^3+11\alpha^2\beta+5\alpha^2\gamma +19\alpha\beta^2
    +24\alpha\beta\gamma +5\alpha\gamma^2+9\beta^3 +19\beta^2\gamma
 +11\beta\gamma^2+\gamma^3 \\
  & & \ \ +\ 11\alpha^2 +39\alpha\beta+28\alpha\gamma +28\beta^2+39\beta\gamma+11\gamma^2 +19\alpha+25\beta+19\gamma+6)n \\
 & & \ \ +\ (\alpha+\beta+\gamma)(\alpha+\beta+\gamma+1)
    (\alpha+\beta+\gamma+2)(\beta+1)] \\
 & & \ \ +\ \left. \rule{0pt}{2.5ex}
  a^2 n(n+\alpha)(\beta+\gamma+2n+1)(\alpha+\beta+\gamma+3n+3)
   \right) \\
  & & +\  \frac{\alpha+\beta+\gamma+2n+1}
  {(\alpha+\beta+\gamma+3n+3)(\alpha+\beta+\gamma+3n+1)^2
   (\alpha+\beta+\gamma+3n)}X_n^{(\alpha,\beta,\gamma)} \\
 & & \ \times\ \left( \rule{0pt}{2.5ex} 12n^3 + (16\alpha+16\beta+10\gamma
 +18)n^2 \right. \\
 & & \ \  +\ [(\alpha+\beta+\gamma)(7\alpha+7\beta+2\gamma)+16\alpha 
 +16\beta+10\gamma]n  \\
 & & \ \ 
+\ (\alpha+\beta+\gamma)^2(\alpha+\beta)+(\alpha+\beta+\gamma)
 (3\alpha+3\beta+2\gamma+2)  \\
 & & \ \ +\ a \left[ 12n^3 + (10\alpha+16\beta+16\gamma
 +18)n^2 \right. \\
 & & \ \ +\ [(\alpha+\beta+\gamma)(2\alpha+7\beta+7\gamma)+10\alpha 
 +16\beta+16\gamma]n  \\
 & & \ \ +\ \left. \rule{0pt}{2.5ex} \left.
   (\alpha+\beta+\gamma)^2(\beta+\gamma)+(\alpha+\beta+\gamma)
 (2\alpha+3\beta+3\gamma+2)\right] \right) \\
 & & -\  \frac{(\alpha+\beta+\gamma+2n+1)^2}
   {(\alpha+\beta+\gamma+3n+1)^2} (X_n^{(\alpha,\beta,\gamma)})^2.
\end{eqnarray*} 
The coefficient $C_{n,n}$ of $x^{2n-3}$ in
$P_{n,n}^{(\alpha,\beta,\gamma)}(x;a)$ can be computed in a similar 
way, and the  coefficient $C_{n+1,n}$ of $x^{2n-2}$ in
$P_{n+1,n}^{(\alpha,\beta,\gamma)}(x;a)$ is given by
\[  C_{n+1,n}^{(\alpha,\beta,\gamma)} =
   C_{n,n}^{(\alpha,\beta+1,\gamma)}-
  X_n^{(\alpha,\beta,\gamma)} B_{n,n}^{(\alpha,\beta,\gamma)}
  \frac{\alpha+\beta+\gamma+2n+1}{\alpha+\beta+\gamma+3n+1}  .  \]
A lengthy but straightforward calculation, using (\ref{eq:d}), then gives
\begin{eqnarray*}
  d_{2n} & = &  \frac{-an(n+\beta)(\alpha+\beta+\gamma+2n)
   (\alpha+\beta+\gamma+2n-1)[n+\gamma + a(n+\alpha)]}
   {(\alpha+\beta+\gamma+3n-2)(\alpha+\beta+\gamma+3n-1)
    (\alpha+\beta+\gamma+3n)^2(\alpha+\beta+\gamma+3n+1)} \\
  & & +\  \frac{n(\alpha+\beta+\gamma+2n)
  (\alpha+\beta+\gamma+2n-1)X_{n-1}^{(\alpha,\beta,\gamma)}}{(\alpha+\beta+\gamma+3n-2)(\alpha+\beta+\gamma+3n-1)
    (\alpha+\beta+\gamma+3n)^2(\alpha+\beta+\gamma+3n+1)} \\
 & & \ \times\ \left[
    (n+\gamma)(\alpha+\beta+2n) -2a(n+\gamma)(n+\alpha) + a^2(n+\alpha)
   (\beta+\gamma+2n) \right] ,
\end{eqnarray*}
and
\begin{eqnarray*}
 d_{2n+1} & = & \frac{n(\alpha+\beta+\gamma+2n+1)(\alpha+\beta+\gamma+2n)}
  {(\alpha+\beta+\gamma+3n+2)(\alpha+\beta+\gamma+3n+1)^2
  (\alpha+\beta+\gamma+3n)^2 (\alpha+\beta+\gamma+3n-1)} \\
  & & \times\ \left( \rule{0pt}{2.5ex}
    (n+\gamma)(\alpha+\beta+2n)(\alpha+\beta+2n+1) \right. \\
  & & \ \ -\ a(n+\alpha)(n+\gamma)(2\alpha+2\beta-\gamma+3n+1) \\
  & & \ \ -\ a^2 (n+\alpha)(n+\gamma)(-\alpha+2\beta+2\gamma+3n+1) \\
  & & \ \ +\ \left. \rule{0pt}{2.5ex}
    a^3(n+\alpha)(\beta+\gamma+2n)(\beta+\gamma+2n+1) \right) \\
 & & -\ \frac{n(\alpha+\beta+\gamma+2n+1)(\alpha+\beta+\gamma+2n)
   X_n^{(\alpha,\beta,\gamma)}}
  {(\alpha+\beta+\gamma+3n+1)^2
  (\alpha+\beta+\gamma+3n)^2 (\alpha+\beta+\gamma+3n-1)} \\
 & & \times\ \left[
(n+\gamma)(\alpha+\beta+2n) -2a (n+\alpha)(n+\gamma)
   +a^2(n+\alpha)(\beta+\gamma+2n) \right]. 
\end{eqnarray*}
The asymptotic behavior of these recurrence coefficients can easily be
found using (\ref{eq:Xinf}), giving
\[   \lim_{n\to \infty} b_{2n} = \frac{a+1}9 + \frac{2x_1}{3}, \quad
     \lim_{n\to \infty} b_{2n+1} = \frac{5(a+1)}9 - \frac{2x_1}{3}, \]
\[ \lim_{n\to\infty} c_{2n} = \frac{4}{81} (a^2-a+1), \quad
   \lim_{n\to\infty} c_{2n+1} = -\frac{4}{9} x_1^2 + \frac{8}{27} x_1
  +\frac{1}{81} (4a^2-a+4), \]
\[ \lim_{n\to\infty} d_{2n} = \frac{4}{243} [2(a^2-a+1)x_1 - a(a+1)], \]
\[  \lim_{n\to\infty} d_{2n+1} = \frac{4}{729} (4a^3-3a^2-3a+4)
      -\frac{8x_1}{243} (a^2-a+1), \]
where $x_1$ is the zero of $g'(x)$ in $[a,0]$ and $g(x)=(x-a)x(x-1)$.
These formulas can be made more symmetric by also using the zero
$x_2$ of $g'(x)$ in $[0,1]$ and using the fact that $x_1+x_2=2(a+1)/3$:
\[ \lim_{n\to \infty} b_{2n} = \frac{a+1}9 + \frac{2x_1}{3}, \quad
   \lim_{n\to \infty} b_{2n+1} = \frac{a+1}9 - \frac{2x_2}{3}, \]
\[     \lim_{n\to\infty} c_{n} = \frac{4}{81} (a^2-a+1), \]
\[  \lim_{n\to\infty} d_{2n} = - \frac{4}{27}  g(x_1), \quad
   \lim_{n\to\infty} d_{2n} = - \frac{4}{27}  g(x_2).   \]

\subsection{Jacobi-Laguerre polynomials}

When we consider the weights $w_1(x) = (x-a)^\alpha |x|^\beta e^{-x}$ on $[a,0]$, with $a < 0$, and
$w_2(x) = (x-a)^\alpha |x|^\beta e^{-x}$ on $[0,\infty)$, then we are again using one weight but on two
touching intervals, one of which is the finite interval $[a,0]$ (Jacobi part), the other the infinite
interval $[0,\infty)$ (Laguerre part). This system was considered by Sorokin \cite{sor1}. 
The corresponding \textit{Jacobi-Laguerre polynomials} $L_{n,m}^{(\alpha,\beta)}(x;a)$
 satisfy the orthogonality relations
\begin{eqnarray*}
  \int_a^0 L_{n,m}^{(\alpha,\beta)}(x;a) (x-a)^\alpha |x|^\beta e^{-x}
   x^k\, dx & = &  0, \qquad k=0,1,\ldots, n-1, \\
  \int_0^\infty L_{n,m}^{(\alpha,\beta)}(x;a) (x-a)^\alpha x^\beta e^{-x}
   x^k\, dx & = &  0, \qquad k=0,1,\ldots, m-1.
\end{eqnarray*}
The raising operator is
\begin{equation}
 \frac{d}{dx} \left[ (x+a)^\alpha x^\beta e^{-x} L_{n,m}^{(\alpha,\beta)}(x;a) \right] 
  = -(x-a)^{\alpha-1} x^{\beta-1} e^{-x} L_{n+1,m+1}^{(\alpha-1,\beta-1)}(x;a),
\end{equation}
from which the Rodrigues formula follows:
\begin{equation}
 \frac{d^m}{dx^m} \left[ (x-a)^{\alpha+m} x^{\beta+m} e^{-x} L_{k,0}^{(\alpha+m,\beta+m)}(x;a) \right] 
  = (-1)^m (x-a)^\alpha x^\beta e^{-x} L_{m+k,m}^{(\alpha,\beta)}(x;a). 
\end{equation}
 From this Rodrigues formula we can proceed as before to find an expression for the polynomials,
but it is more convenient to view these Jacobi-Laguerre polynomials as a limit case of the
Jacobi-Angelesco polynomials  
\begin{equation}
  L_{n,m}^{(\alpha,\beta)}(x;a) = \lim_{\gamma \to \infty}
   \gamma^{n+m} P_{n,m}^{(\alpha,\beta,\gamma)}(x/\gamma;a/\gamma),
\end{equation}
so that (\ref{eq:jacang2}) gives
\begin{equation}  \label{jaclag2}
  L_{n,n}^{(\alpha,\beta)}(x;a) = \sum_{k=0}^n \sum_{j=0}^{n-k}
        \binom{\alpha+n}{k} \binom{\beta+n}{j} \frac{(-1)^{k+j}(x-a)^{n-k} x^{n-j}}{(n-k-j)!}.
\end{equation}
For the recurrence coefficients in
\[  xP_n(x) = P_{n+1}(x) + b_n P_n(x) + c_n P_{n-1}(x) + d_n P_{n-2}(x), \]
where $P_{2n}(x) = L_{n,n}^{(\alpha,\beta)}(x;a)$ and $P_{2n+1}(x) = L_{n+1,n}^{(\alpha,\beta)}(x;a)$
we have in terms of the corresponding recurrence coefficients of the Jacobi-Angelesco polynomials
\begin{eqnarray*}
  b_n & = & \lim_{\gamma\to \infty} \gamma    b_n^{(\alpha,\beta,\gamma)}(a/\gamma), \\
  c_n & = & \lim_{\gamma\to \infty} \gamma^2          c_n^{(\alpha,\beta,\gamma)}(a/\gamma), \\
   d_n & = & \lim_{\gamma\to \infty} \gamma^3             d_n^{(\alpha,\beta,\gamma)}(a/\gamma),
\end{eqnarray*}
and
\[  \lim_{\gamma \to \infty} \gamma X_n^{(\alpha,\beta,\gamma)}(a/\gamma)
   = \frac{ \int_a^0 x(x-a)^{\alpha+n} |x|^{\beta+n} e^{-x} \, dx}
   {\int_a^0 (x-a)^{\alpha+n} |x|^{\beta+n} e^{-x} \, dx} :=
 X_n^{(\alpha,\beta)} . \]
This gives
\begin{eqnarray*}
   b_{2n} & = & n + X_n^{(\alpha,\beta)}, \\
  b_{2n+1} & = & 3n+\alpha+\beta+2+a-X_n^{(\alpha,\beta)}, \\
  c_{2n} & = & n(\alpha+\beta+2n), \\
  c_{2n+1} & = & n(\alpha+\beta+2n+1)-a(n+\beta+1) + (\alpha+\beta+2n+2+a) X_n^{(\alpha,\beta)} - (X_n^{(\alpha,\beta)})^2 , \\
  d_{2n} & = & -an(\beta+n) + n(\alpha+\beta+2n) X_{n-1}^{(\alpha,\beta)}, \\
  d_{2n+1} & = & n[(\alpha+\beta+2n)(\alpha+\beta+2n+1)+a(n+\alpha)]
   -n(\alpha+\beta+2n) X_n^{(\alpha,\beta)}. 
\end{eqnarray*}
For large $n$ we have $X_n^{(\alpha,\beta)} =a/2 + o(1)$ so that
\begin{eqnarray*}
   \lim_{n \to \infty} \frac{b_{n}}{n} & = & \begin{cases} 1/2 & \mathrm{if\ } n \equiv 0 \pmod{2}, \\
                                                       3/2 & \mathrm{if\ } n \equiv 1 \pmod{2},
                                         \end{cases} \\
  \lim_{n \to \infty} \frac{c_n}{n^2} & = & 1/2, \\
  \lim_{n \to \infty} \frac{d_n}{n^3} & = & \begin{cases} 0 & \mathrm{if\ } n \equiv 0 \pmod{2}, \\
                                                       1/2 & \mathrm{if\ } n \equiv 1 \pmod{2}.
                                         \end{cases}
\end{eqnarray*}

\subsection{Laguerre-Hermite polynomials}
Another limit case of the Jacobi-Angelesco polynomials are the multiple orthogonal polynomials
$H_{n,m}^{(\beta)}(x)$ for which
\begin{eqnarray*}
  \int_{-\infty}^0 H_{n,m}^{(\beta)}(x) |x|^\beta e^{-x^2}
   x^k\, dx & = &  0, \qquad k=0,1,\ldots, n-1, \\
  \int_0^\infty H_{n,m}^{(\beta)}(x)  x^\beta e^{-x^2}
   x^k\, dx & = &  0, \qquad k=0,1,\ldots, m-1.
\end{eqnarray*}
We call these \textit{Laguerre-Hermite polynomials} because both weights are supported
on semi-infinite intervals (Laguerre) with a common weight that resembles
the Hermite weight. These polynomials were already considered (for general $r$) by
Sorokin \cite{sor3}. The limit case is obtained by taking
\begin{equation}
  H_{n,m}^{(\beta)}(x) = \lim_{\alpha \to \infty}
   (\sqrt{\alpha})^{n+m} P_{n,m}^{(\alpha,\beta,\alpha)}
  (x/\sqrt{\alpha};-1).
\end{equation}
This allows us to obtain the raising operator, the Rodrigues formula, an explicit expression, and
the recurrence coefficients by taking the appropriate limit passage in the formulas for the
Jacobi-Angelesco polynomials. For the recurrence coefficients this gives
\begin{eqnarray*}
  b_n & = & \lim_{\alpha\to \infty} \sqrt{\alpha}    b_n^{(\alpha,\beta,\alpha)}(a=-1), \\
  c_n & = & \lim_{\alpha\to \infty} \alpha          c_n^{(\alpha,\beta,\alpha)}(a=-1), \\
   d_n & = & \lim_{\alpha\to \infty} (\sqrt{\alpha})^3             d_n^{(\alpha,\beta,\alpha)}(a=-1),
\end{eqnarray*}
and
\[  \lim_{\alpha \to \infty} \sqrt{\alpha} X_n^{(\alpha,\beta,\alpha)}(a=-1)
   = \frac{ \int_{-\infty}^0 x |x|^{\beta+n} e^{-x^2} \, dx}
   {\int_{-\infty}^0  |x|^{\beta+n} e^{-x^2} \, dx} :=
 X_n^{(\beta)} , \]
from which we find
\begin{eqnarray*}
   b_{2n} & = & X_n^{(\beta)}, \\
  b_{2n+1} & = & -X_n^{(\beta)}, \\
  c_{2n} & = & n/2, \\
  c_{2n+1} & = & \frac{2n+\beta+1}{2}-(X_n^{(\beta)})^2 , \\
  d_{2n} & = & \frac{n}{2} X_{n-1}^{(\beta)}, \\
  d_{2n+1} & = & \frac{-n}{2} X_n^{(\beta)}. 
\end{eqnarray*}
For large $n$ we have
\[    X_n^{(\beta)} =-\sqrt{\frac{\beta+n}2} + o(\sqrt{n}), \]
so that
\begin{eqnarray*}
   \lim_{n \to \infty} \frac{b_{n}}{\sqrt{n}} & = & \begin{cases} -1/2 & \mathrm{if\ } n \equiv 0 \pmod{2}, \\
                                                       1/2 & \mathrm{if\ } n \equiv 1 \pmod{2},
                                         \end{cases} \\
  \lim_{n \to \infty} \frac{c_n}{n} & = & 1/4, \\
  \lim_{n \to \infty} \frac{d_n}{(\sqrt{n})^3} & = & \begin{cases} -1/8 & \mathrm{if\ } n \equiv 0 \pmod{2}, \\
                                                       1/8 & \mathrm{if\ } n \equiv 1 \pmod{2}.
                                         \end{cases}
\end{eqnarray*}

\section{Open research problems}

In the previous sections we gave a short description of multiple orthogonal polynomials
and a few  examples. For a more detailed account of multiple orthogonal
polynomials we refer to Aptekarev \cite{apt} and Chapter 4 of the book of Nikishin
and Sorokin \cite{nikisor}. Multiple orthogonal polynomials arise naturally
in Hermite-Pad\'e approximation of a system of (Markov) functions. For this kind
of simultaneous rational approximation we refer to Mahler \cite{mahler} and
de Bruin \cite{bruin1}, \cite{bruin2}. Hermite-Pad\'e approximation goes back to the nineteenth century,
and many algebraic aspects have been investigated since then: existence and
uniqueness, recurrences, normality of indices, etc. The more detailed analytic investigation
of the zero distribution, the $n$th root asymptotics, and the strong asymptotics
is more recent and mostly done by researchers from the schools around Nikishin 
\cite{niki1}, \cite{niki2} and Gonchar \cite{gonrak}, \cite{gonraksor}. See in particular the work
of Aptekarev \cite{apt}, Kalyagin \cite{kalya1}, \cite{kalyaron}, 
Bustamante and L\'opez \cite{lopez}, but also the
work by Driver and Stahl \cite{drivst1}, \cite{drivst2} and Nuttall \cite{nut}. First one needs to understand the analysis
of ordinary orthogonal polynomials, and then one has a good basis for studying
this extension, for which there are quite a few possibilities for research.

\subsection{Special functions}
The research of orthogonal polynomials as special functions has now
 led to a classification
and arrangement of various important (basic hypergeometric) orthogonal polynomials. In Section
3 we gave a few multiple orthogonal polynomials of the same flavor as the very
classical orthogonal polynomials of Jacobi, Laguerre, and Hermite. Regarding these
very classical multiple orthogonal polynomials, a few open problems arise:

\begin{enumerate}
\item Are the polynomials given in Section 3 the only possible very classical multiple
orthogonal polynomials? The answer very likely is no. First one needs to make clear
what the notion of classical multiple orthogonal polynomial means. A possible way
is to start from a Pearson type equation for the weights. If one chooses one weight
but restricted to disjoint intervals, as we did for the Jacobi-Angelesco, Jacobi-Laguerre, 
and Laguerre-Hermite polynomials, then Aptekarev, Marcell\'an and Rocha \cite{aptmarroc}
used the Pearson equation for this weight as the starting point of their characterization.
For several weights it is more natural to study a Pearson equation for the vector
of weights $(w_1,w_2,\ldots,w_r)$. Douak and Maroni \cite{doumar1},
\cite{doumar2} have given a complete characterization of all type II multiple orthogonal
polynomials for which the derivatives are again type II multiple
orthogonal polynomials (Hahn's characterization for the Jacobi,
Laguerre, and Hermite polynomials, and the Bessel polynomials
if one allows moment functionals which are not positive definite). 
They call such polynomials  classical $d$-orthogonal polynomials,
where $d$ corresponds to our $r$, i.e., the number of
weights (functionals) needed for the orthogonality.
Douak and Maroni show that this class of multiple orthogonal
polynomials is characterized by a Pearson equation of the form
\[  (\Phi \vec{w} )' + \Psi \vec{w} = \vec{0}, \]
where $\vec{w} = (w_1,\ldots,w_r)^t$ is the vector of weights,
and $\Psi$ and $\Phi$ are $r\times r$ matrix polynomials:
\[  \Psi(x) = \begin{pmatrix}
     0 & 1 & 0 & \cdots & 0 \\
     0 & 0 & 2 & \cdots & 0 \\
     \vdots & \vdots & \vdots & \ddots & \vdots \\
     0 & 0 & 0 & \cdots & r-1 \\
     \psi(x) & c_1 & c_2 & \cdots & c_{r-1}
   \end{pmatrix} , \]
with $\psi(x)$ a polynomial of degree one and $c_1,\ldots,c_{r-1}$ constants, and 
\[  \Phi(x) = \begin{pmatrix}
     \phi_{1,1}(x) & \phi_{1,2}(x) & \cdots & \phi_{1,r}(x) \\
     \phi_{2,1}(x) & \phi_{2,2}(x) & \cdots & \phi_{2,r}(x) \\
       \vdots & \vdots & \cdots & \vdots \\
     \phi_{r,1}(x) & \phi_{r,2}(x) & \cdots & \phi_{r,r}(x) 
   \end{pmatrix} , \]
where $\phi_{i,j}(x)$ are polynomials of degree at most two. In fact
only $\phi_{r,1}$ can have degree at most two and
all other polynomials are constant or of degree one, depending
on their position in the matrix $\Phi$. 
Douak and Maroni actually investigate the more general case where
orthogonality is given by $r$ linear functionals, rather than by
$r$ positive measures.
We believe that Hahn's characterization is not the appropriate
property to define classical multiple orthogonal polynomials, but gives
a more restricted class. None of the seven families, given in the
present paper, belong to the class studied by Douak and Maroni, but
their class certainly contains several interesting families of multiple orthogonal polynomials. In fact, the matrix Pearson equation could 
result from a single weight (and its derivatives) satisfying a
higher order differential equation with polynomial coefficients. As an example, one can have
multiple orthogonal polynomials with weights
$w_1(x) = 2 x^{\alpha+\nu/2} K_\nu(2\sqrt{x})$ and
$w_2(x) = 2 x^{\alpha+(\nu+1)/2} K_{\nu+1}(2\sqrt{x})$ on $[0,\infty)$, where
$K_\nu(x)$ is a modified Bessel function and $\alpha > -1$, $\nu \geq 0$ (see \cite{wvayak} and \cite{cheik}). 
\item The polynomials of Jacobi, Laguerre, and Hermite all satisfy a linear second order differential equation of Sturm-Liouville type. A possible
way to extend this characterizing property is to look for multiple
orthogonal polynomials satisfying a linear differential
equation of order $r+1$. Do the seven families in this paper
have such a differential equation? If the answer is yes,  
then an explicit construction would be desirable. We only
worked out in detail the case where $r=2$, so the search is for a third order
differential equation for all the polynomials considered in Section 3.
 Such a third order equation has been found for certain
Jacobi-Angelesco systems in \cite{kalyaron}. For the Angelesco systems in Section 3 this
third order differential equation indeed exists and it was constructed in
\cite{aptmarroc}. The existence (and construction) is open for the AT systems.
A deeper problem is to characterize all the multiple orthogonal polynomials satisfying
a third order (order $r+1$) differential equation, extending Bochner's result
for ordinary orthogonal polynomials. Observe that we already know appropriate raising
operators for the seven systems described in Section 3. If one can construct lowering
operators as well, then a combination of the raising and lowering operators will give
the differential equation, which will immediately be in factored form. Just
differentiating will usually not be sufficient (except for the class
studied by Douak and Maroni): if we take $P_{n,m}'(x)$, then this is a polynomial
of degree $n+m-1$, so one can write it as $P_{n-1,m}(x)+$ lower order terms, but  also
as $P_{n,m-1}(x)+$ lower order terms. So it is not clear which of the multi-indices
has to be lowered. Furthermore, the lower order terms will not vanish in general since there
usually are not enough orthogonality conditions to make them disappear.
\item In the present paper we only considered the type II multiple orthogonal polynomials.
Derive explicit expressions and relevant properties of the corresponding vector 
$(A_{n,m}(x), B_{n,m}(x))$ of type
I multiple orthogonal polynomials. Type I and type II multiple orthogonal polynomials are
connected by
 \[    P_{n,m}(x) = \mathrm{const.} \begin{vmatrix}  A_{n+1,m}(x) & B_{n+1,m}(x) \\
                                     A_{n,m+1}(x) & B_{n,m+1}(x)
                    \end{vmatrix}, \]
but from this it is not so easy to obtain the type I polynomials.
\item So far we limited ourselves to the very classical orthogonal polynomials of
Jacobi, Laguerre, and Hermite. Discrete orthogonal polynomials, such as those of Charlier, Kravchuk, Meixner, and Hahn, can also be considered and several
kinds of discrete multiple orthogonal polynomials can be worked out. It would not be a good
idea to do this case by case, since these polynomials are all connected by limit 
transitions, with the Hahn polynomials as the starting family. At a later stage, one
could also consider multiple orthogonal polynomials on a quadratic lattice and
on the general exponential lattice, leading to $q$-polynomials. Again, all these
families are related, with the Askey-Wilson polynomials as the family from which
all others can be obtained by limit transitions. Do these polynomials
have a representation as a (basic) hypergeometric function? Recall that
we needed an Appell hypergeometric function of two variables for
the Jacobi-Angelesco polynomials, so that one may need to consider
(basic) hypergeometric functions of several variables.  
\item Multiple orthogonal polynomials arise naturally in the study
of Hermite-Pad\'e approximation, which is simultaneous rational
approximation to a vector of $r$ functions. In this respect it is quite
natural to study multiple orthogonal polynomials as orthogonal vector polynomials. This approach is very useful in trying to extend
results for the case $r=1$ to the case $r>1$ by looking for an
appropriate formulation using vector algebra. Van Iseghem 
already used this approach to formula a vector QD-algorithm
for multiple orthogonal polynomials \cite{iseg}. Several algebraic
aspects of multiple orthogonal polynomials follow easily from
the vector orthogonality \cite{soriseg1}, \cite{labeck}. A further generalization is to
study matrix orthogonality, where the matrix need not be a square matrix
\cite{soriseg2}. Orthogonal polynomials and Pad\'e approximants
are closely related to certain continued fractions (J-fractions
and S-fractions). For multiple orthogonal polynomials there is
a similar relation with vector continued fractions and the
Jacobi-Perron algorithm \cite{par}. The seven families which we considered in this
paper lead to seven families of vector continued fractions, which could
be studied in more detail in the framework of continued fractions.
Finally, one may wonder whether it is possible to use hypergeometric functions
of matrix argument in the study of multiple orthogonal polynomials? 
\end{enumerate}

\subsection{Non-symmetric banded operators}
In Section 2 the connection between multiple orthogonal polynomials and banded Hessenberg operators
of the form
\[  \setcounter{MaxMatrixCols}{11}
\begin{pmatrix} 
  a_{0,0} & 1    \\
  a_{1,1} & a_{1,0} & 1 \\
  a_{2,2} & a_{2,1} & a_{2,0} & 1 \\
   \vdots & & & \ddots &\ddots \\
  a_{r,r} & a_{r,r-1} & \cdots & & a_{r,0} & 1 \\
          & a_{r+1,r} & \ddots &  & & a_{r+1,0} & 1 \\
          &   & \ddots & \ddots & &  & \ddots & \ddots \\
           &  &   &  \ddots & \ddots & & & \ddots & 1 \\
          &   &   &     & a_{n,r} & a_{n,r-1} & \cdots & a_{n,1} &  a_{n,0} & \ddots  \\
          &   &   &     &          & \ddots  & \cdots & & & \cdots & \ddots
 \end{pmatrix} \]
was explained. For ordinary orthogonal polynomials the operator is tridiagonal and can always
be made symmetric, and often it can be extended in a unique way to a self-adjoint operator
(e.g, when all the coefficients are bounded). The spectrum of this tridiagonal operator corresponds
to the support of the orthogonality measure, and the spectral measure is precisely the
orthogonality measure. Each tridiagonal matrix with ones on the upper diagonal and positive 
coefficients on the lower diagonal, corresponds to a system of orthogonal polynomials on the
real line (Favard's theorem). 
Some preliminary work on the spectral theory of the higher order  operators ($r >1$) was done
by Kalyagin \cite{kalya2}, \cite{kalya3}, \cite{kalya4}, \cite{aptkal}, but there are still
quite a few open problems here.

\begin{enumerate}
\item What is the proper extension of Favard's theorem for these higher order banded Hessenberg
operators? Not every banded Hessenberg operator corresponds to a system of multiple
orthogonal polynomials with orthogonality relations on the real line. There needs to be
additional structure, but so far this additional structure is still unknown. There is a weak
version of the Favard theorem that gives multiple orthogonality with respect to linear
functionals (\cite{iseg}, \cite{kalya5}), but a stronger version that gives positive measures on the real line
is needed. How do we recognize an Angelesco system, an AT system, or one of the combinations
considered in \cite{gonraksor} from the recurrence coefficients (from the operator)? The special case where all the diagonals are zero, except for the
upper diagonal (which contains 1's) and the lower diagonal, has been
studied in detail in \cite{aptkalvis}. They show that when the
lower diagonal contains positive coefficients, the operator corresponds
to multiple orthogonal polynomials on an $(r+1)$-star in the complex
plane. Using a symmetry transformation, similar to the quadratic
transformation that transforms Hermite polynomials to Laguerre polynomials,
this also gives an AT system of multiple orthogonal polynomials on $[0,\infty)$.
\item The asymptotic behavior of the recurrence coefficients of the seven
 systems described above
is known. Each of the limiting operators deserves to be investigated in more detail. The
limiting operator for the Jacobi-Pi\~neiro polynomials is a Toeplitz operator, and hence can
be investigated in more detail. See, e.g., \cite{wva2} for this case. 
Some of the other limiting
operators are block Toeplitz matrices and can be investigated as well.
Are there any multiple orthogonal polynomials having such recurrence coefficients? The
Chebyshev polynomials of the second kind have this property when one deals with tridiagonal
operators.
\item The next step would be to work out a perturbation theory, where one allows certain
perturbations of the limiting matrices. Compact perturbations would be the first step, trace
class perturbations would allow us to give more detailed results.
\end{enumerate} 

\subsection{Applications} 
\begin{enumerate}
\item Hermite-Pad\'e approximation was introduced by Hermite for his proof
of the transcendence of $e$. More recently it became clear that
Ap\'ery's proof of the irrationality of $\zeta(3)$ relies 
on an AT system of multiple orthogonal polynomials with weights
$w_1(x) = 1$, $w_2(x) = -\log(x)$ and $w_3(x) = \log^2(x)$ on 
$[0,1]$. These multiple orthogonal polynomials are basically
limiting cases of Jacobi-Pi\~neiro polynomials where $\alpha_0=0=\alpha_1=\alpha_2$. 
A very interesting problem is to prove irrationality of other
remarkable constants, such as $\zeta(5)$, Catalan's constant, or
Euler's constant. Transcendence proofs will even be better.
See \cite{apt} \cite{wva1} for the connection between multiple orthogonal
polynomials, irrationality, and transcendence.
\item In numerical analysis one uses
orthogonal polynomials when one constructs Gauss quadrature. In a similar
way one can use multiple orthogonal polynomials to construct optimal
quadrature formulas for jointly approximating $r$ integrals of the
same function $f$ with respect to $r$ weights $w_1,\ldots,w_r$.
See, e.g., Borges \cite{borges}, who apparently is not aware that he is
using multiple orthogonal polynomials. Gautschi \cite{gaut} has
summarized some algorithms for computing  recurrence
coefficients, quadrature nodes (zeros of orthogonal polynomials)
and quadrature weights (Christoffel numbers) for ordinary Gauss
quadrature. A nice problem is to modify these algorithms so that they
compute recurrence coefficients, zeros of multiple orthogonal
polynomials (eigenvalues of banded Hessenberg operators) and
quadrature weights for simultaneous Gauss quadrature.
\end{enumerate}

\end{document}